\documentclass[11pt,a4paper]{article}
\usepackage{fullpage}

\usepackage[utf8x]{inputenc} 
\usepackage{amsfonts} 
\usepackage{amsmath,amsthm,amssymb} 

\usepackage{cite} 

\usepackage{graphicx} 

\usepackage{hyperref}

\usepackage{color}

\def\q{\boldsymbol{q}}

\def\sgn{\text{sgn}}
\def\Id{\boldsymbol{I}}

\def\eps{\varepsilon}

\newcommand{\vr}{\varrho}

\newcommand{\vrs}{\varrho^*}
\newcommand{\pt}{\partial_{t}}

\newcommand{\Div}{{\nabla\,\cdot\, }}

\newcommand{\vu}{\boldsymbol{v}}

\newcommand{\vc}[1]{{\boldsymbol #1}}

\newcommand{\Grad}{\nabla}
\newcommand{\lr}[1]{\left( #1 \right)}
\newcommand{\ep}{\varepsilon}

\newcommand{\eq}[1]{\begin{equation}
\begin{split}
#1
\end{split}
\end{equation}}
\newcommand{\eqh}[1]{\begin{equation*}
\begin{split}
#1
\end{split}
\end{equation*}}

\newtheorem{thm}{Theorem}

\newtheorem{prop}[thm]{Proposition}

\title{Finite Volume approximations of the Euler system\\ with variable congestion}

\author{Pierre Degond$^1$, Piotr Minakowski$^2$, Laurent Navoret$^3$, Ewelina Zatorska$^{1}$}
\date{}

\begin{document}

\maketitle
\centerline{1. Department of Mathematics, Imperial College London, }
\centerline{London SW7 2AZ, United Kingdom.}
\bigskip
  \centerline{2.  Interdisciplinary Center for Scientific Computing, Heidelberg University,}
 \centerline{Im Neunheimer Feld 205, D-69120 Heidelberg, Germany.}
\bigskip

\centerline{3. Institut de Recherche Math\'ematique Avanc\'ee, UMR 7501,}
\centerline{ Universit\'e de Strasbourg et CNRS, 7 rue Ren\'e Descartes,}
\centerline{ 67000 Strasbourg, France \& Inria TONUS}

\bigskip

 
\begin{abstract} We are interested in the numerical simulations of the Euler system with variable congestion encoded by a singular pressure \cite{DeMiZa2016}. This model describes for instance the macroscopic motion of a crowd with individual congestion preferences. We propose an asymptotic preserving (AP) scheme based on a conservative formulation of the system in terms of density, momentum and density fraction. A second order accuracy version of the scheme is also presented. We validate the scheme on one-dimensionnal test-cases and compare it with a scheme previously proposed in \cite{DeMiZa2016} and extended here to higher order accuracy. We finally carry out two dimensional numerical simulations and show that the model exhibit typical crowd dynamics.

 
\end{abstract}

\section{Introduction}
In this work we study two phase compressible/incompressible Euler system with variable congestion:
\begin{subequations}\label{sysSl}
\begin{align}
 \pt\vr+ \Div (\vr\vu) = 0,\label{rho} \\
 \partial_t (\vr\vu) + \Div (\vr\vu \otimes \vu) + \Grad\pi +\Grad p\lr{\frac{\vr}{\vr^*}} = \vc{0},\label{mom}\\
 \pt\vr^*+\vu\cdot\Grad\vr^*=0,\label{rho_star}\\
 0\leq \vr\leq\vrs,\label{cons0}\\
 \pi(\vrs-\vr)=0,\quad \pi\geq 0,\label{pineq0}
\end{align}
 \end{subequations}
 with the initial data
\eq{\label{initial_rho}
\vr(0,x)=\vr_0(x)\geq 0,\quad\vu(0,x)=\vu_0(x),\quad \vrs(0,x)=\vrs_0(x),\quad \vr_0<\vrs_0,}
where the unknowns are: $\vr=\vr(t,x)$ -- the mass density, $\vu=\vu(t,x)$ -- the velocity, $\vrs=\vrs(t,x)$ -- the congestion density, and $\pi$ -- the congestion pressure. The barotropic pressure $p$ is an explicit function of the density fraction $\frac{\vr}{\vr^*}$
\eq{\label{pbar}
p\left(\frac{\vr}{\vr^*}\right)=\lr{\frac{\vr}{\vr^*}}^\gamma,\quad \gamma>1,
}
and plays the role of the background pressure. 

The congestion pressure $\pi$ appears only when the density $\vr$ achieves its maximal value, the congestion density $\vrs$. Therefore $\vrs$ is sometimes referred to as the barrier or the threshold density. It was observed in \cite{LM99}, and then generalized in \cite{DeMiZa2016}, that the restriction on the density \eqref{cons0} is equivalent with the condition 
\eq{
  \Div\vu =0 \ \text{in} \ \{\vr=\vrs\},\label{div0}
}
if only $\vr, \vu,\vrs$ are sufficiently regular solutions of the continuity equation \eqref{rho} and the transport equation \eqref{rho_star}. For that reason, system \eqref{sysSl} can be seen as a free boundary problem for the interface between the compressible (uncongested) regime $\{\vr<\vrs\}$ and the incompressible (congested) regime $\{\vr=\vrs\}$.

The main purpose of this work is to analyze \eqref{sysSl} numerically, i.e. to propose the numerical scheme capturing the phase transition. To this end we use the fact that \eqref{sysSl} can be obtained as a limit when $\ep\to 0$ of the compressible Euler system with the congestion pressure $\pi$  replaced by a singular approximation $\pi_\ep$ 
\eq{\label{pie}
\pi_\ep\lr{\frac{\vr}{\vr^*}}=\ep \lr{\frac{\frac{\vr}{\vr^*}}{1-\frac{\vr}{\vr^*}}}^\alpha, \quad \alpha >0.
}
Note that for fixed $\ep>0$, $\pi_\ep\to\infty$ when $\vr\to \vrs$. Therefore, at least formally, for $\ep\to 0$, $\pi_\ep$ converges to a measure supported on the set of singularity, i.e. $\{(x,t)\in \Omega\times (0,T):\vr(x,t)=\vrs(x,t)\}$. The rigorous proof of this fact is an open problem, at least for the Euler type of systems. There have been, however, several results for a viscous version of the model, see \cite{BrPeZa} for the one-dimensional case, \cite{PeZa} for multi-dimensional domains and space-dependent congestion $\vrs(x)$ and \cite{DeMiZa2016} for the case of congestion density satisfying the transport equation \eqref{rho_star}. An alternative approximation leading to a similar  two-phase system was considered first by P.-L. Lions and N. Masmoudi \cite{LM99}, and more recently for the model of tumour growth \cite{PeVa15}. The advantage of approximation \eqref{pie} considered here lies in the fact that for each $\ep$ fixed, the solutions to the approximate system stay in the physical regime, i.e. $\vr\leq\vrs$. This feature is especially important for the numerical purposes, see for example \cite{Ma12} for further discussion on this subject.

System \eqref{sysSl} is a generalization of the pressureless Euler system with the maximal density constraint
\begin{subequations}\label{sys0}
\begin{align}
 \pt\vr+ \Div (\vr\vu) = 0, \\
 \partial_t (\vr\vu) + \Div (\vr\vu \otimes \vu) + \Grad\pi  = \vc{0},\\
 0\leq \vr\leq 1\\
 \pi(\vr-1)=0,\quad \pi\geq 0.
\end{align}
 \end{subequations}
introduced originally by Bouchut et al. \cite{BoBrCoRi}, who also proposed the first numerical scheme based on an approach developed earlier for the pressureless systems, see for example \cite{Brenier84}, and the projection argument. The model was studied later on by Berthelin \cite{Be02, Berthelin16} by passing to the limit in the so-called  sticky-blocks dynamics, see also \cite{Wolansky}, and a very interesting recent paper \cite{PeWe2016} using the Lagrangian approach for the monotone rearrangement of the solution to prove the existence of solutions to \eqref{sys0} with additional memory effects. 

The pressureless Euler equations with the density constraint were originally introduced in order to describe the motion of particles of finite size. Our model extends this concept by including the variance of the size of particles. In system  \eqref{sysSl} $\vrs$ is given initially and is transported along with the flow. 

One can also think of $\vrs$ as a congestion preference of individuals moving in the crowd (cars, pedestrians), which is one of the factors determining their final trajectory and the speed of motion. The macroscopic modelling of crowd is one of possible approaches and it allows to determine the averaged quantities such as the density and the mean velocity rather than the precise position of an individual. One of the first models of this kind based on classical mechanics was introduced by Henderson \cite{henderson1971}. More sophisticated  model was introduced by Hughes \cite{Hughes2002} where the author considers the continuity equation equipped with a  phenomenological constitutive relation between the velocity and the density.  For a survey of the crowd models we refer the reader to  \cite{Colombo2005, Bellomo2008, Tosin2009, Maury2010, Jiang2010} and to the review paper \cite{bellomo2011}. 

As far as the numerical methods are concerned, the macroscopic models of pedestrian flow with condition preventing the overcrowding were studied, for example in \cite{TwGoDu13}. The influence of the maximal density constraint  was investigated also in the context of vehicular traffic in \cite{Berthelin_M3AS_08, BeBr}. The strategy that we want to adapt in this paper, i.e. to use the singularities of the pressure similar to \eqref{pie} has been developed in the past for a number of Euler-like systems for the traffic models \cite{BeBr, Berthelin_M3AS_08, Berthelin_ARMA_08}, collective dynamics \cite{DeHuNa, DeHu}, or granular flow \cite{Maury07, Perrin16}. 
In our previous work  \cite{DeMiZa2016}, we have drafted the numerical scheme for system \eqref{sysSl} in the one-dimensional case. We used a splitting algorithm at each time step that consists of three sub-steps. At first, the hyperbolic part is solved with the AP-preserving method presented in \cite{DeHuNa}. Next the diffusion is solved by means of cell-centered finite volume scheme, and the transport of the congested density is resolved with the upwind scheme.

The extension of this method to two-dimensions is one of the main results of the present paper. We also propose an alternative scheme using different formulation in terms of the {\emph{conservative variables}}: the density $\vr$, the momentum $\q = \vr\vu$, and the density fraction $Z=\frac{\vr}{\vrs}$: 
\begin{subequations} \label{sysSZ}
 \eq{\label{sysSZa}
 \pt\vr+ \Div \, \q = 0,}
  \eq{\label{sysSZb}\partial_t \q + \Div \lr{\frac{\q\otimes \q}{\vr} + \pi_\ep(Z)\Id + p(Z)\Id}  = \vc{0},}
  \eq{\label{sysSZc}\pt Z+\Div\lr{Z\frac{\q}{\vr}}=0,}
with the initial data
\eq{\label{initialZ}
\vr(0,x)=\vr_0(x),\quad\q(0,x)=\vc{q}_0(x),\quad Z(0,x)=Z_0(x),}
\end{subequations}
where $Z _0= \frac{\vr_0}{\vrs_0}$, and $\vc{q}_0=\vr_0\vu_0$. $\Id$ denotes the identity tensor. This is a stricly hyperbolic system whose wave speeds in the $x_{1}$-direction are given by:
\eq{
&\lambda_{1}^{\eps}(\vr,q_{1}, Z) =  \frac{q_{1}}{\vr} - \sqrt{\frac{Z}{\vr} p_{\eps}'(Z)},\, \\
&\lambda_{2}^{\eps}(\vr,q_{1}, Z) = \frac{q_{1}}{\vr}, \\
&\lambda_{3}^{\eps}(\vr,q_{1}, Z) = \frac{q_{1}}{\vr}+ \sqrt{\frac{Z}{\vr} p_{\eps}'(Z)}, \label{eq:eigenval}
}
where  $p_\ep=p+\pi_\ep$, and $q_{1}$ denotes the component of $\q$ in the $x_{1}$ direction. Consequently, in region where the density $\vr$ is closely congested, i.e. $Z$ is close to 1,  the characteristic speeds of the system are extremely large. This corresponds to the nearly incompressible dynamics.

The paper is organised as follows. In Section \ref{sec:numscheme} we present our numerical schemes using the two formulations \eqref{sysSl} and \eqref{sysSZ}. They are referred to as  $(\vr, \q)$-method/SL and $(\vr, \q, Z)$-method, respectively. In Section~\ref{sec:z1st} we describe 
 the first-order semi-discretization in time and the full discretization for the $(\vr, \q, Z)$-method.  Then, in Section~\ref{sec:second_order_discretization}, we discuss the second order scheme for the $(\vr, \q, Z)$-method. At last,  in Section~\ref{sec:semiLagscheme} we present the $(\vr, \q)$-method/SL for the system written in terms of the physical variables \eqref{sysSl}. Section \ref{sec:validation} is devoted to validation of the schemes on the Riemann problem whose solutions are described in \ref{sec:solRiemann}. Finally, in Section~\ref{sec:2d_comp} we discuss the two-dimensional numerical results: in Section \ref{sec:test} we present how these schemes work for three different initial congestion densities, and in Section~\ref{sec:application} we present an application of $(\vr, \q)$-method/SL  to model crowd behaviour in the evacuation scenario.

\section{Numerical schemes}\label{sec:numscheme}

In this section, we first introduce a numerical scheme based on system \eqref{sysSZ} using the conservative variables. In order to use large time steps not restricted by too drastic CFL condition, implicit-explicit (IMEX) type methods need to be designed. The scheme can be solved through the following steps: first an elliptic equation on the density fraction $Z$ is solved, and then we update $\q$ and $\vr$, respectively. 

Such scheme is compared with an extension of the method introduced in \cite{DeMiZa2016}, where the congestion density is advected separately from the update of $\vr$ and $\q$. For the sake of completeness, a description of the scheme is given in Section \ref{sec:semiLagscheme}.

\subsection{The first order $(\vr, \q, Z)$-method}\label{sec:z1st}

\paragraph{Discretization in time}
We adopt the previous work  \cite{DeHuNa} to introduce a method treating implicitly the stiff congestion pressure $\pi_{\eps}(Z)$. We consider a constant time step $\Delta t > 0$ and $\vr^{n}$, $\q^{n}$, $Z^{n}$, $\vr^{\ast\, n}$ denote the approximate solution at time $t^{n} = n \Delta t$, $\forall n \in \mathbb{N}$. We thus consider the following semi-implicit time discretization:
\begin{subequations}\label{eq:semidiscrete}
\begin{align}
&\frac{\vr^{n+1} - \vr^{n}}{\Delta t} + \nabla_{x} \cdot \q^{n+1} = 0,\label{eq:semidiscrete_rho}\\
&\frac{\q^{n+1} - \q^{n}}{\Delta t} + \nabla_{x}\cdot\left(\frac{\q^{n}\otimes \q^{n}}{\vr^{n}} + p(Z^{n}) \Id \right) + \nabla_{x}(\pi_\ep(Z^{n+1})) = 0,\label{eq:semidiscrete_q}\\
&\frac{Z^{n+1} - Z^{n}}{\Delta t} + \nabla_{x}\cdot \left(Z^{n} \frac{\q^{n+1}}{\vr^{n}}\right) = 0\label{eq:semidiscrete_z}.
\end{align}
\end{subequations}
Note that in the flux term in equation \eqref{eq:semidiscrete_z}, the momentum is taken implicitly. Inserting \eqref{eq:semidiscrete_q} into \eqref{eq:semidiscrete_z}, we obtain:
\begin{align*}
&\frac{Z^{n+1} - Z^{n}}{\Delta t} + \nabla_{x}\cdot \left(Z^{n} \frac{\q^{n}}{\vr^{n}}\right) \\
&\qquad - \Delta t\, \nabla_{x}\cdot \left(\frac{Z^{n}}{\vr^{n}}\nabla_{x}\cdot\left(\frac{\q^{n}\otimes \q^{n}}{\vr^{n}}+ p(Z^{n}) \Id \right) + \frac{Z^{n}}{\vr^{n}} \nabla_{x}(\pi_\ep(Z^{n+1}))\right) = 0,
\end{align*}
This is an elliptic equation on the unknown $Z^{n+1}$, that can be written as:
\begin{equation}
Z^{n+1} - \Delta t^{2}\, \nabla_{x}\cdot \left(\frac{Z^{n}}{\vr^{n}}\nabla_{x}\left(\pi_\ep(Z^{n+1})\right)\right) = \phi(\vr^n,\
q^n,Z^n), \label{eq:ellipticZsemidiscrete}
\end{equation}
where
\begin{align*}
&\phi(\vr^n,\
q^n,Z^n) \\
&= Z^{n} + \Delta t^{2}\, \nabla_{x}\cdot \left(\frac{Z^{n}}{\vr^{n}}\nabla_{x}\cdot\left(\frac{\q^{n}\otimes \q^{n}}{\vr^{n}} + p(Z^{n})\Id\right)\right) - \Delta t\, \nabla_{x}\cdot \left(Z^{n} \frac{\q^{n}}{\vr^{n}}\right).
\end{align*}
The $n$-th time step of the scheme is decomposed into three parts: first get $Z^{n+1}$ when solving \eqref{eq:ellipticZsemidiscrete}, then compute $\q^{n+1}$ thanks to \eqref{eq:semidiscrete_q} and then $\vr^{n+1}$ from \eqref{eq:semidiscrete_rho}.

\paragraph{Discretization in space}
We only derive the fully discrete scheme in the one-dimensional case; the two-dimensional formula are given in \ref{sec:scheme_2d}. We consider the computational domain $[0,1]$ and a spatial space step $\Delta x = 1/N_{x} > 0$, with $N_{x} \in \mathbb{N}$: the mesh points are thus $x_{i} = i \Delta x$, $\forall i \in \left\{0,\ldots, N_{x}\right\}$. Let $\vr^{n}_{i}$, $\q^{n}_{i}$, $Z^{n}_{i}$, $\vr^{\ast\, n}_{i}$ denote the approximate solution at time $t^{n}$ on mesh cell $[x_{i},x_{i+1}]$. The spatial discretization have to capture correctly the entropic solutions of the hyperbolic system. To derive the fully discrete scheme, we thus make the same algebra on the following fully discrete system:
\begin{subequations}
\begin{align}
&\frac{\vr^{n+1}_{i} - \vr^{n}_{i}}{\Delta t} + \frac{1}{\Delta x} (F_{i+\frac{1}{2}}^{n+1} - F_{i-\frac{1}{2}}^{n+1})  = 0,\label{eq:discrete_rho}\\
&\frac{q^{n+1}_{i} - q^{n}_{i}}{\Delta t} + \frac{1}{\Delta x} (G_{i+\frac{1}{2}}^{n} - G_{i-\frac{1}{2}}^{n})  + \frac{\pi_\ep(Z^{n+1}_{i+1}) - \pi_\ep(Z^{n+1}_{i-1})}{2\Delta x} = 0,\label{eq:discrete_q}\\
&\frac{Z^{n+1}_{i} - Z^{n}_{i}}{\Delta t} + \frac{1}{\Delta x} (H_{i+\frac{1}{2}}^{n+1} - H_{i-\frac{1}{2}}^{n+1})  = 0\label{eq:discrete_z}.
\end{align}
\label{eq:fulldiscrete}
\end{subequations}
where the stiff pressure is discretized by the centered finite difference and the numerical fluxes $F^{n+1}$, $G^{n}$, $H^{n+1}$ (we denote implicit-explicit fluxes by current timestep $n+1$ and fully explicit fluxes by previous timestep $n$) are splitted into centered part and the upwinded part:  
\begin{align}
&F_{i+\frac{1}{2}}^{n+1} = \frac{1}{2}\big(q^{n+1}_{i+1}+q^{n+1}_{i}\big) - (D_{\vr})_{i+\frac{1}{2}}^{n},\label{eq:flux_rho}\\
&G_{i+\frac{1}{2}}^{n} = \frac{1}{2}\Big(\frac{(q_{i+1}^{n})^{2}}{\vr_{i+1}^{n}}+\frac{(q_{i}^{n})^{2}}{\vr_{i}^{n}} + p(Z_{i+1}^{n})+ p(Z_{i}^{n})\Big)  - (D_{q})_{i+\frac{1}{2}}^{n},\label{eq:flux_q}\\
&H_{i+\frac{1}{2}}^{n+1} = \frac{1}{2}\Big( \frac{Z^{n}_{i+1}}{\vr^{n}_{i+1}} q^{n+1}_{i+1}+ \frac{Z^{n}_{i}}{\vr^{n}_{i}} q^{n+1}_{i}\Big)  - (D_{Z})_{i+\frac{1}{2}}^{n}.\label{eq:flux_Z}
\end{align}
The upwinded parts are given explicitly. They can be given by the diagonal Rusanov (or local Lax-Friedrichs) upwindings:
\begin{equation}
(D_{w})_{i+\frac{1}{2}}^{n} = \frac{1}{2} c^n_{i+\frac{1}{2}} \big(w_{i+1}^{n} - w_{i}^{n}\big),\label{eq:upwindings}
\end{equation}
for any conserved quantities $w$, where $c^n_{i+\frac{1}{2}}$ is the maximal characteristic speed (in absolute value):
\begin{equation}
c^n_{i+\frac{1}{2}} = \max \left\{ \left|\lambda_{k}^{0}\big(\vr^n_{i+1}, q^n_{i+1}, Z^n_{i+1}\big)\right|, \left|\lambda_{k}^{0}\big(\vr^n_{i}, q^n_{i}, Z^n_{i}\big)\right|,\quad k = 1, 2, 3\right\},
\label{eq:maxcharspeed}
\end{equation}
where $\lambda_{k}^{0}$ are given by eq. \eqref{eq:eigenval} with $\eps = 0$ (no congestion pressure). These correspond to the eigenvalues of the hyperbolic system taken explicitly in \eqref{eq:semidiscrete}. One could also consider less diffusive numerical fluxes like the Polynomial upwind scheme \cite{DePeyRuVill}.

Like in the semi-discrete case, we now obtain the fully discrete elliptic equation on $Z$ by replacing the implicit momentum terms appearing in the flux $H$ \eqref{eq:flux_Z} by their expressions given by the momentum equation \eqref{eq:discrete_q}. We get:
\begin{align*}
&Z^{n+1}_{i} - Z^{n}_{i} + \frac{\Delta t}{\Delta x} (\bar H_{i+1/2}^{n} - \bar H_{i-1/2}^{n})  \\
&- \frac{\Delta t^{2}}{\Delta x^{2}} \frac{1}{2}\Big( \frac{Z^{n}_{i+1}}{\vr^{n}_{i+1}} (G_{i+\frac{3}{2}}^{n} - G_{i+\frac{1}{2}}^{n})- \frac{Z^{n}_{i-1,j}}{\vr^{n}_{i-1}} (G_{i-\frac{1}{2}}^{n} - G_{i-\frac{3}{2}}^{n})\Big)\\
&-  \frac{\Delta t^{2}}{\Delta x^{2}} \frac{1}{2}\Big( \frac{Z^{n}_{i+1}}{\vr^{n}_{i+1}} \big(\pi_\ep(Z^{n+1}_{i+2}) - \pi_\ep(Z^{n+1}_{i})\big) - \frac{Z^{n}_{i-1}}{\vr^{n}_{i-1}} \big(\pi_\ep(Z^{n+1}_{i}) - \pi_\ep(Z^{n+1}_{i-2})\big)\Big) = 0,
\end{align*}
where $\bar H^{n}$ denotes the same expression as \eqref{eq:flux_Z} where all quantities are taken explicitly:
\begin{equation*}
\bar H_{i+\frac{1}{2}}^{n} = \frac{1}{2}\Big( \frac{Z^{n}_{i+1}}{\vr^{n}_{i+1}} q^{n}_{i+1}+ \frac{Z^{n}_{i}}{\vr^{n}_{i}} q^{n}_{i}\Big)  - (D_{Z})_{i+\frac{1}{2}}^{n}.
\end{equation*}  
In practice, in order to preserve the constraint $Z \leqslant 1$,  this elliptic equation is solved with respect to the congestion pressure variable $\pi_{\eps}$: 
\begin{align}
Z_i^{n+1}((\pi_{\ep})^{n+1}_{i}) -  \frac{\Delta t^{2}}{\Delta x^{2}} \frac{1}{2} &\Big( \frac{Z^{n}_{i+1}}{\vr^{n}_{i+1}} \big[(\pi_\eps)^{n+1}_{i+2} - (\pi_\ep)^{n+1}_{i}\big] \label{eq:Z}\\
& - \frac{Z^{n}_{i-1}}{\vr^{n}_{i-1}} \big[(\pi_\ep)^{n+1}_{i} - (\pi_\ep)^{n+1}_{i-2}\big]\Big) = \phi(\vr^{n}, q^{n}, Z^{n})_{i},\nonumber
\end{align}
where the right-hand side is given by:
\begin{align}
\phi(\vr^{n}, q^{n}, Z^{n})_{i} = &Z^{n}_{i}  - \frac{\Delta t}{\Delta x} (H_{i+1/2}^{n} - H_{i-1/2}^{n}) \label{eq:Z2}  \\
&+ \frac{\Delta t^{2}}{\Delta x^{2}} \frac{1}{2}\Big( \frac{Z^{n}_{i+1}}{\vr^{n}_{i+1}} (G_{i+\frac{3}{2}}^{n} - G_{i+\frac{1}{2}}^{n})- \frac{Z^{n}_{i-1,j}}{\vr^{n}_{i-1}} (G_{i-\frac{1}{2}}^{n} - G_{i-\frac{3}{2}}^{n})\Big).\nonumber
\end{align}
and $Z(\pi_{\ep})$ is the inverse function of $\pi_{\ep}(Z)$. This equation is supplemented by periodic or Dirichlet boundary conditions and the non-linear equation is solved using the Newton iterations.

The $(n\!+\!1)$-th time step of the algorithm thus consists in getting $Z^{n+1}$ by solving \eqref{eq:Z}-\eqref{eq:Z2} and then obtaining $q^{n+1}$ from \eqref{eq:discrete_q} and $\vr^{n+1}$ from \eqref{eq:discrete_rho}.

Since the singular pressure $\pi_{\eps}$ is treated implicitly, the scheme remains stable even for small $\varepsilon$. The stability condition only depends on the wave speeds of the explicit part of the scheme, that is under the Courant-Friedrichs-Levy (CFL) condition:
\begin{equation}
\Delta t \leqslant \frac{\Delta x}{\underset{j = 1, 2, 3;\, x \in [0,1], t \in [0,T]}{\max}\left\{ |\lambda_{j}^{0}(x,t)| \right\}},
\label{eq:CFL}
\end{equation}
where $\lambda_{j}^{0}$, given by eq. \eqref{eq:eigenval}, denotes the eigenvalues of the hyperbolic system with no congestion pressure ($\eps = 0$). The scheme is asymptotically stable with respect to $\varepsilon$.

\subsection{The second order $(\vr, \q, Z)$-method}\label{sec:second_order_discretization}
\paragraph{Discretization in time}
The second-order discretization in time is based on the combined Runge-Kutta 2 / Crank-Nicolson (RK2CN) method as described in \cite{CorDeKum}: it consists into replacing Euler explicit by Runge-Kutta 2 solver and Euler Implicit by Crank-Nicolson solver in semi-discretization \eqref{eq:semidiscrete}. We here only detail the semi-discretized scheme. However, to be unambiguous, we will denote by $\mathcal{D}_{\vr}$, $\mathcal{D}_{q}$, and $\mathcal{D}_{Z}$ the numerical diffusion terms resulting from the upwinding terms and the divergence operators will be replaced by centered fluxes. We thus consider the following scheme:

\noindent{\it{First step}} (half time step): get $\vr^{n+1/2}$, $\q^{n+1/2}$ and $Z^{n+1/2}$ from
\begin{subequations}
\begin{align}
&\frac{\vr^{n+1/2} - \vr^{n}}{\Delta t/2} + \nabla_{x} \cdot \q^{n+1/2} - \mathcal{D}_{\vr}^{n} = 0,\\
&\frac{\q^{n+1/2} - \q^{n}}{\Delta t/2}\!+\! \nabla_{x}\cdot\left(\frac{\q^{n}\otimes \q^{n}}{\vr^{n}} + p(Z^{n}) \Id \right) \!-\! \mathcal{D}_{q}^{n} \!+\! \nabla_{x}(\pi_\ep(Z^{n+1/2})) \!=\! 0,\label{eq:semidiscrete_q2}\\
&\frac{Z^{n+1/2} - Z^{n}}{\Delta t/2} + \nabla_{x}\cdot \left(\frac{Z^{n}}{\vr^{n}} \q^{n+1/2}\right) - \mathcal{D}_{Z}^{n}= 0\label{eq:semidiscrete_z2}.
\end{align}
\end{subequations}  
{\it{Second step}} (full time step): get $\vr^{n+1}$, $\q^{n+1}$ and $Z^{n+1}$ from
\begin{subequations}
\begin{align}
&\frac{\vr^{n+1} - \vr^{n}}{\Delta t} + \nabla_{x} \cdot \left(\frac{\q^{n+1}+\q^{n}}{2}\right) - \mathcal{D}_{\vr}^{n} = 0,\\
&\frac{\q^{n+1} - \q^{n}}{\Delta t} + \nabla_{x}\cdot\left(\frac{\q^{n+1/2}\otimes \q^{n+1/2}}{\vr^{n+1/2}} + p(Z^{n+1/2}) \Id \right) -  \mathcal{D}_{q}^{n+1/2} \nonumber\\
&\hspace{5cm}+ \nabla_{x}\left(\frac{\pi_\ep(Z^{n}) +\pi_\ep(Z^{n+1})}{2}\right) = 0,\label{eq:semidiscrete_q3}\\
&\frac{Z^{n+1} - Z^{n}}{\Delta t} + \nabla_{x}\cdot \left(\frac{Z^{n+1/2}}{\vr^{n+1/2}} \frac{\q^{n+1}+\q^{n}}{2}\right) - \mathcal{D}_{Z}^{n} = 0\label{eq:semidiscrete_z3}.
\end{align}
\end{subequations}  
Like in the first-oder scheme, equations \eqref{eq:semidiscrete_q2}-\eqref{eq:semidiscrete_z2} and \eqref{eq:semidiscrete_q3}-\eqref{eq:semidiscrete_z3} result in elliptic equations for $Z$. In practice, the scheme may fail capturing discontinuities, in particular when small values of $\eps$ are concerned. Indeed,  the semi-implicit pressure $\big(\pi_\ep(Z^{n}) +\pi_\ep(Z^{n+1})\big)/2$ in \eqref{eq:semidiscrete_q3} is constrained to be larger than $\pi_\ep(Z^{n})/2$ preventing from having large discontinuities in pressure. One way to overcome this difficulty is to dynamically replace this semi-implicit pressure by an implicit pressure $\pi_\ep(Z^{n+1})$ as soon as the non-linear solver of the elliptic equation detects a pressure lower than half the explicit one.

\paragraph{Discretization in space}
To get second order accuracy in space, we consider a MUSCL strategy. For any conserved quantity $v$, it consists in introducing at each mesh interface left and right values $w_{L}$ and $w_{R}$:
\begin{align*}
&w_{i, L} = v_{i} + \frac{1}{2}\, \text{minmod}(w_{i}-w_{i-1}, w_{i+1}-w_{i}),\\
&w_{i,R} = v_{i} - \frac{1}{2}\, \text{minmod}(w_{i}-w_{i-1}, w_{i+1}-w_{i}),
\end{align*}
where the minmod function is defined as: 
$$\text{minmod}(a,b) = 0.5\, (\sgn(a)+\sgn(b)) \min(|a|, |b|) .$$ 
Then all explicit terms in fluxes \eqref{eq:flux_rho}-\eqref{eq:flux_q}-\eqref{eq:flux_Z} depend 
on $(\vr_{i,R}^{n}, q_{i,R}^{n}, Z_{i,R}^{n})$ and $(\vr_{i+1, L}^{n}, q_{i+1,L}^{n}, Z_{i+1,L}^{n})$ instead of  $(\vr_{i}^{n}, q_{i}^{n}, Z_{i}^{n})$ and $(\vr_{i+1}^{n}, q_{i+1}^{n}, Z_{i+1}^{n})$. Implicit terms are unchanged in order to be able to get the elliptic equation.

\subsection{Congested Euler/Semi-Lagrangian scheme ($(\vr, \vc{q})$-method/SL)}
\label{sec:semiLagscheme}
\paragraph{Discretization in time}
We consider a scheme based on the non-conservative form \eqref{sysSl} of the congestion transport. This idea was proposed in \cite{DeHuNa} in the context of constant congestion and in \cite{DeMiZa2016} in the context of variable congestion.
 The time-discretization reads:
\begin{subequations}
\begin{align}
&\frac{\vr^{n+1} - \vr^{n}}{\Delta t} + \nabla_{x} \cdot \q^{n+1} = 0,\label{eq:semidiscrete2_rho}\\
&\frac{\q^{n+1} - \q^{n}}{\Delta t} + \nabla_{x}\cdot\left(\frac{\q^{n}\otimes \q^{n}}{\vr^{n}} + p\left(\frac{\vr^{n}}{{\vr^*}^n}\right) \Id \right) + \nabla_{x}\pi_\ep\left(\frac{\vr^{n+1}}{{\vr^*}^n}\right) = 0,\label{eq:semidiscrete2_q}\\
&\frac{{\vr^*}^{n+1} - {\vr^*}^{n}}{\Delta t} + \frac{\q^{n+1}}{\vr^{n+1}}\cdot\nabla_{x} {\vr^*}^{n} = 0\label{eq:semidiscrete2_z}.
\end{align}
\label{eq:semidiscrete2}
\end{subequations}
Inserting\eqref{eq:semidiscrete2_q} into \eqref{eq:semidiscrete2_rho} results in
\begin{align}\label{eq:ellipticRHOsemidiscrete}
\vr^{n+1}  - &\Delta t^{2}\, \Delta_{x}\big(\pi_\ep(\vr^{n+1}/\vr^{\ast\, n})\big)= \nonumber\\
& \vr^{n} - \Delta t \nabla_{x} \cdot \q^{n} + \Delta t^{2}\,   \nabla_{x} \cdot \nabla_{x}\cdot\left(\frac{\q^{n}\otimes \q^{n}}{\vr^{n}} + p(\vr^{n}/\vr^{\ast\, n}) \Id \right).
\end{align}
This is an elliptic equation on the density $\vr^{n+1}$. The $n$-th time step of the scheme is decomposed into three parts: first get $\vr^{n+1}$ when solving \eqref{eq:ellipticRHOsemidiscrete}, then compute $\q^{n+1}$ thanks to \eqref{eq:semidiscrete2_q} and then ${\vr^*}^{n+1}$ from \eqref{eq:semidiscrete2_z}.

\paragraph{Discretization in space} Like for the previous schemes, we restrict the description to the one-dimensional case.
Finite volume discretization is used for the spatial discretization of \eqref{eq:semidiscrete2_rho}-\eqref{eq:semidiscrete2_q} as in section \ref{sec:z1st}, see also \cite{DeHuNa}.  A semi-Lagrangian method is used to solve \eqref{eq:semidiscrete2_z} and thus update the congestion density $\vr^{\ast}$. The congestion density $\vr^{\ast n+1}_{i}$ at node $x_{i}$ and time $t^{n+1}$ is computed as follows: first we integrate back the characteristic line over one time step and then we interpolate the maximal density $\vr^{\ast n}$ at that point. Using Euler scheme for the first step, we obtain:
\begin{equation*}
\vr_{i}^{\ast n+1} = \left[\Pi \vr^{\ast n}\right](x_{i} - q_{i}/\vr_{i}\, \Delta t)
\end{equation*}
where $\Pi \vr^{\ast n}$ is an interpolation function built from the points $(x_{i}, \vr^{\ast n}_{i})$. We here perform a Lagrange interpolation on the $2r + 2$ neighboring points:
$$[\Pi \vr^{\ast}]_{|[x_{i},x_{i+1}]} = \Pi_{\text{Lagrange}}\Big((x_{j},\vr^{\ast}_{j}),\quad i-r+1\leq j\leq i+r\Big).$$
resulting in $2r+1$-th spatial accuracy. First $(r=0)$ and third $(r=1)$ order in space semi-Lagrangian scheme will be used. For more details, we refer to \cite{FalconeFerretti}.

\paragraph{The second order scheme}
Extension of the full scheme to second order accuracy in space is made using the MUSCL strategy for the finite volume fluxes. Extension to second order accuracy in time requires a Crank-Nicolson/Runge Kutta 2 method for $(\vr, q)$ and a second order in time integration of the characteric line for the semi-Lagrangian scheme (with for instance Taylor expansion) combined to a Strang splitting, see \ref{sec:secondorder_rhoscheme}.

\section{One dimensional validation of the schemes}\label{sec:validation}

\subsection{Riemann test-case}

We compare the numerical schemes on one-dimensional Riemann test-cases: the initial data is a discontinuity between two constant states and the solutions are given by the superposition of waves separating constant states. In \ref{sec:solRiemann}, we give the form of these solutions with respect to the relative position of left and right states in the phase space. In the case of colliding states, explicit solutions can be numerically obtained. We thus consider the following Riemann test-case:
\begin{equation}
(\vr_{0}(x), q_{0}(x), \vr^{\ast}_{0}(x)) = \begin{cases}
(\vr_{\ell}, q_{\ell}, \vr^{\ast}_{\ell}) =  (0.7,0.8,1.2), &\text{if } x \leqslant 0.5,\\
 (\vr_{r}, q_{r}, \vr^{\ast}_{r}) =  (0.7, - 0.8,1), &\text{if } x > 0.5.
 \end{cases}
 \label{eq:riempb}
 \end{equation}
 on the domain $[0,1]$. The solution is made of two shock waves and an intermediate contact wave, see \eqref{eq:solutionRiemPb}. The CFL condition \eqref{eq:CFL} can be estimated by:
\begin{equation*}
\Delta t \leqslant \frac{\Delta x}{\underset{ x \in [0,1], t \in [0,T]}{\max} |v(x,t)| + \sqrt{\gamma/ \min \vr^{\ast}(x,t)}}.
\end{equation*}
For the current Riemann test-case with $\gamma = 2$, the time step should satisfy $\Delta t \leqslant 0.4 \Delta x$.

\begin{figure}
\begin{center}
\hspace{-0.6cm}\includegraphics[width=0.45\textwidth]{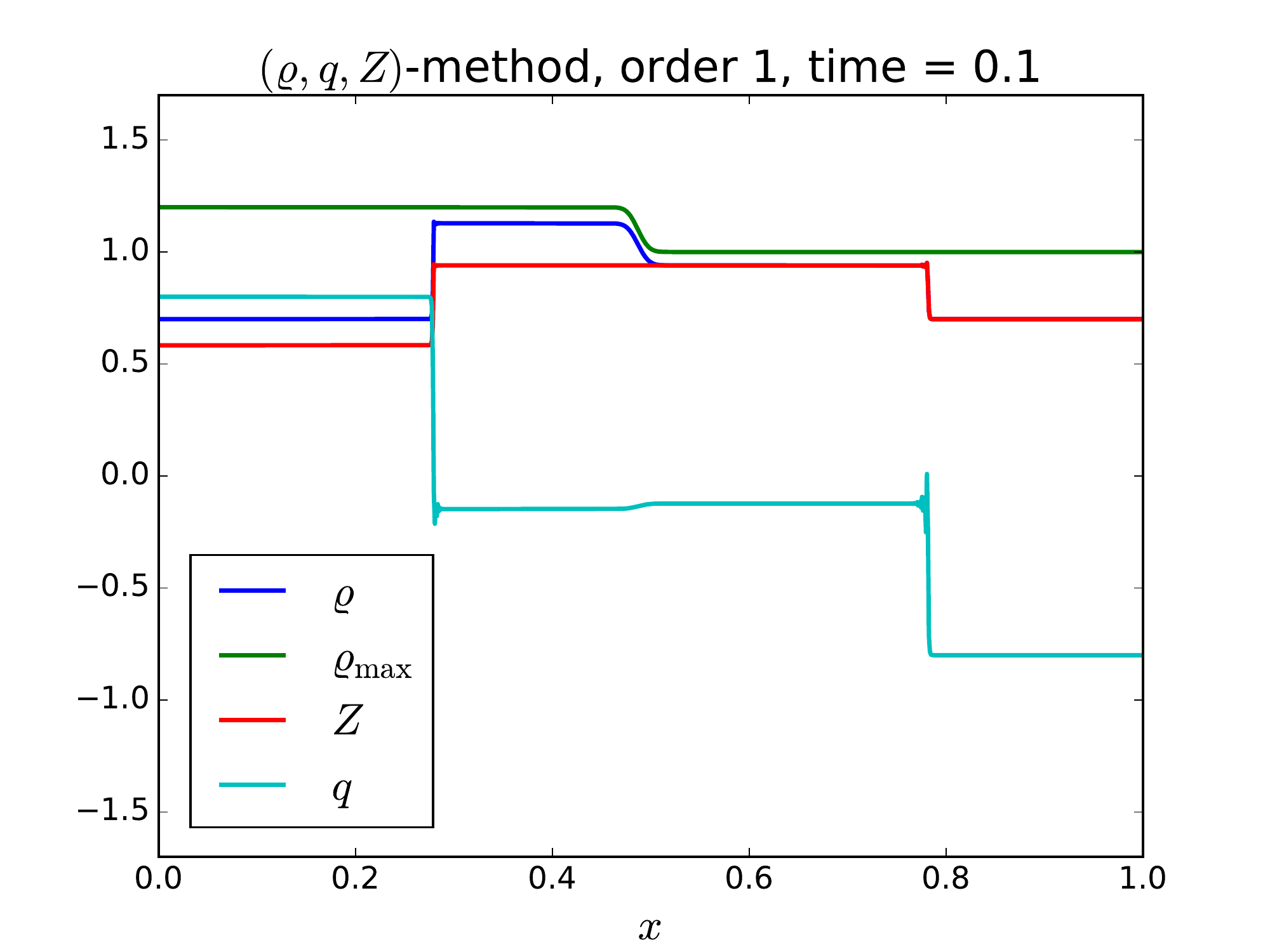}\includegraphics[width=0.45\textwidth]{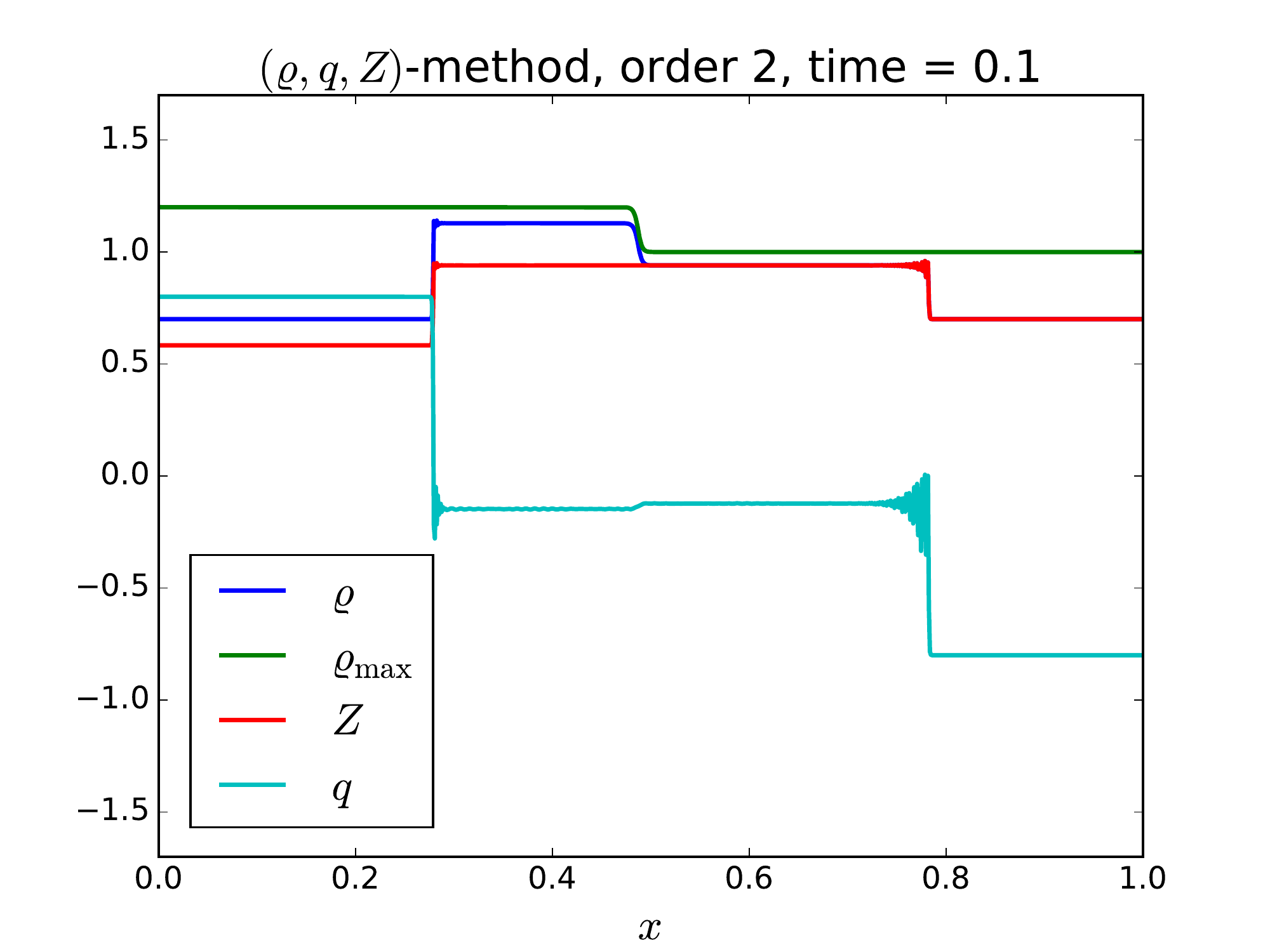}

\hspace{-0.6cm}\includegraphics[width=0.45\textwidth]{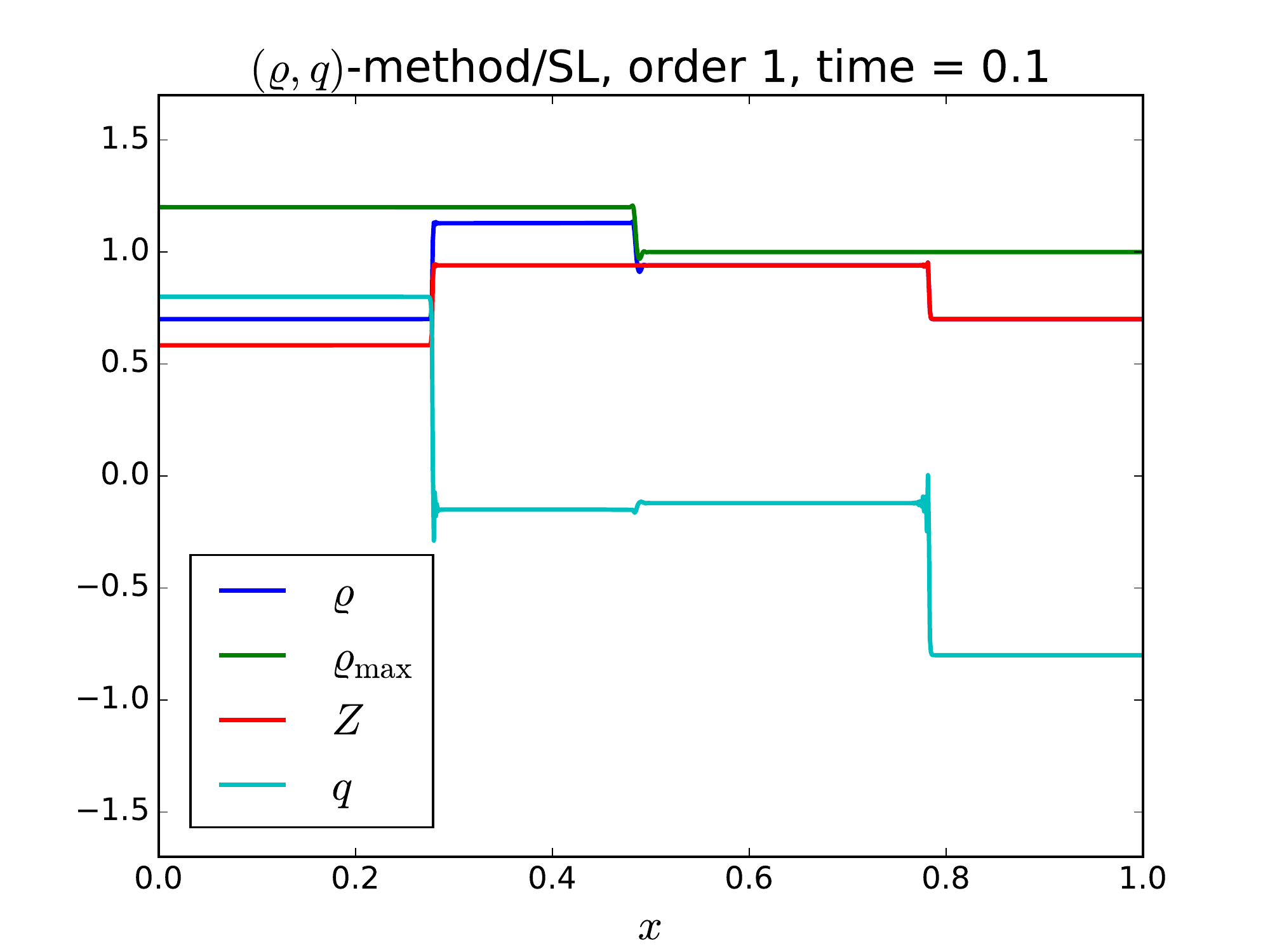}\includegraphics[width=0.45\textwidth]{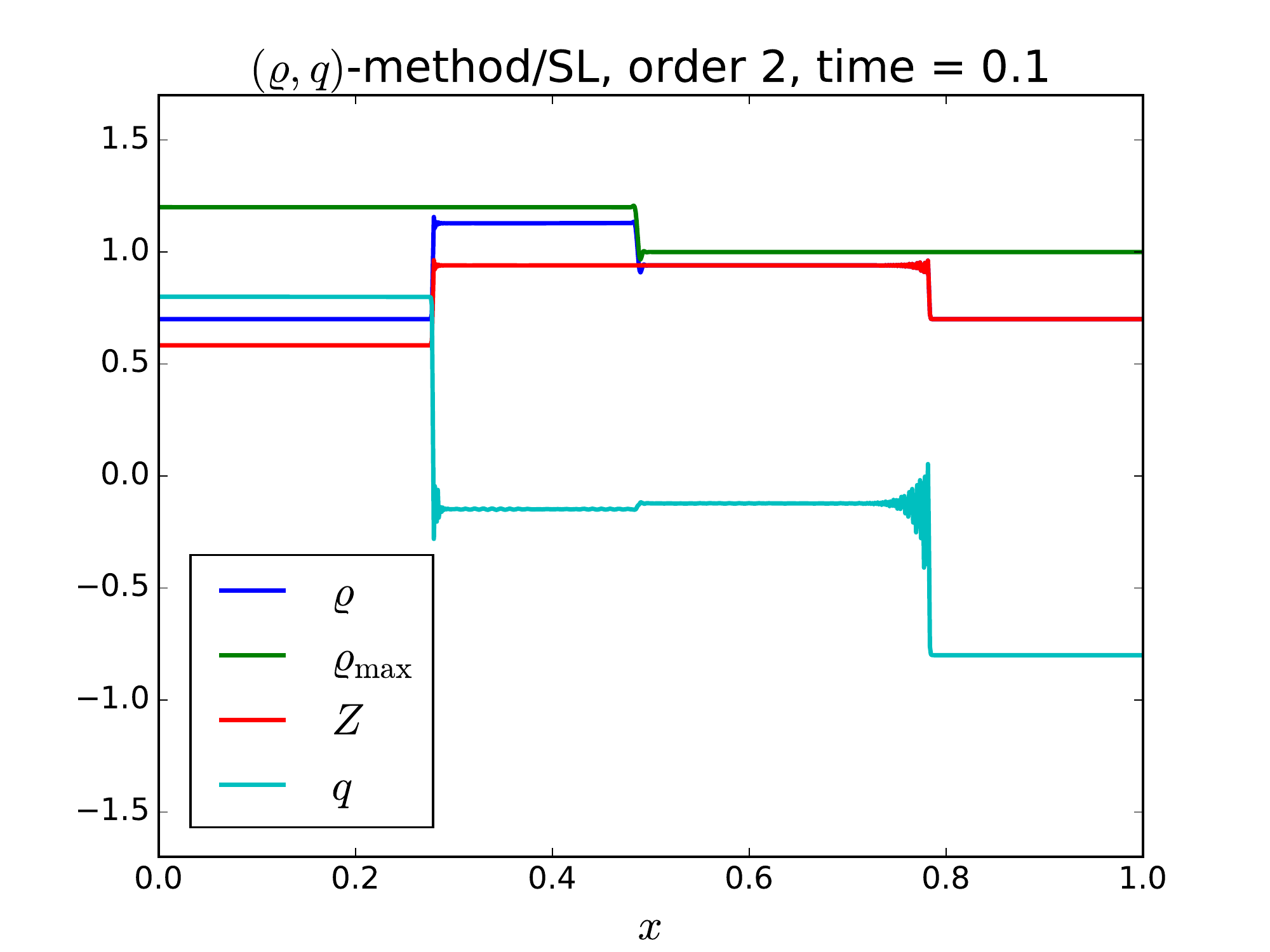}
\end{center}
\caption{Approximate solution to Riemann problem \eqref{eq:riempb} at time $t = 0.1$. Numerical parameters: $\Delta x= 1\times 10^{-3}$, $ \Delta t = 0.1\, \Delta x$, $\alpha = 2$, $\gamma = 2$, $\eps=10^{-2}$.}
\label{sec:riemtestcase}
\end{figure}

\paragraph{Comparison of the schemes ($\eps = 10^{-2}$)}  In Figure~\ref{sec:riemtestcase}, we represent the solution at time $t=0.1$ with the different schemes using $\Delta t = 0.1 \Delta x$. The $(\vr, \q, Z)$-method refers to the method introduced in Section \ref{sec:z1st} for the first order and in Section \ref{sec:second_order_discretization} for the second order scheme. The $(\vr, q)$-method/SL refers to the method described in Section \ref{sec:semiLagscheme}. For the latter scheme, we use the third order semi-Lagrangian scheme for the transport of the congestion density $\vr^\ast$. 

We observe that all the methods correctly capture the exact solution. The $(\vr,\q)$-method/SL better captures the contact discontinuity at $x \approx 0.487$ since we use a third order accurate scheme for the transport of $\vr^{\ast}$. Limiters could be used to avoid overshoot and undershoot at this location.

\begin{figure}
\begin{center}
\hspace{-0.6cm}\includegraphics[width=0.45\textwidth]{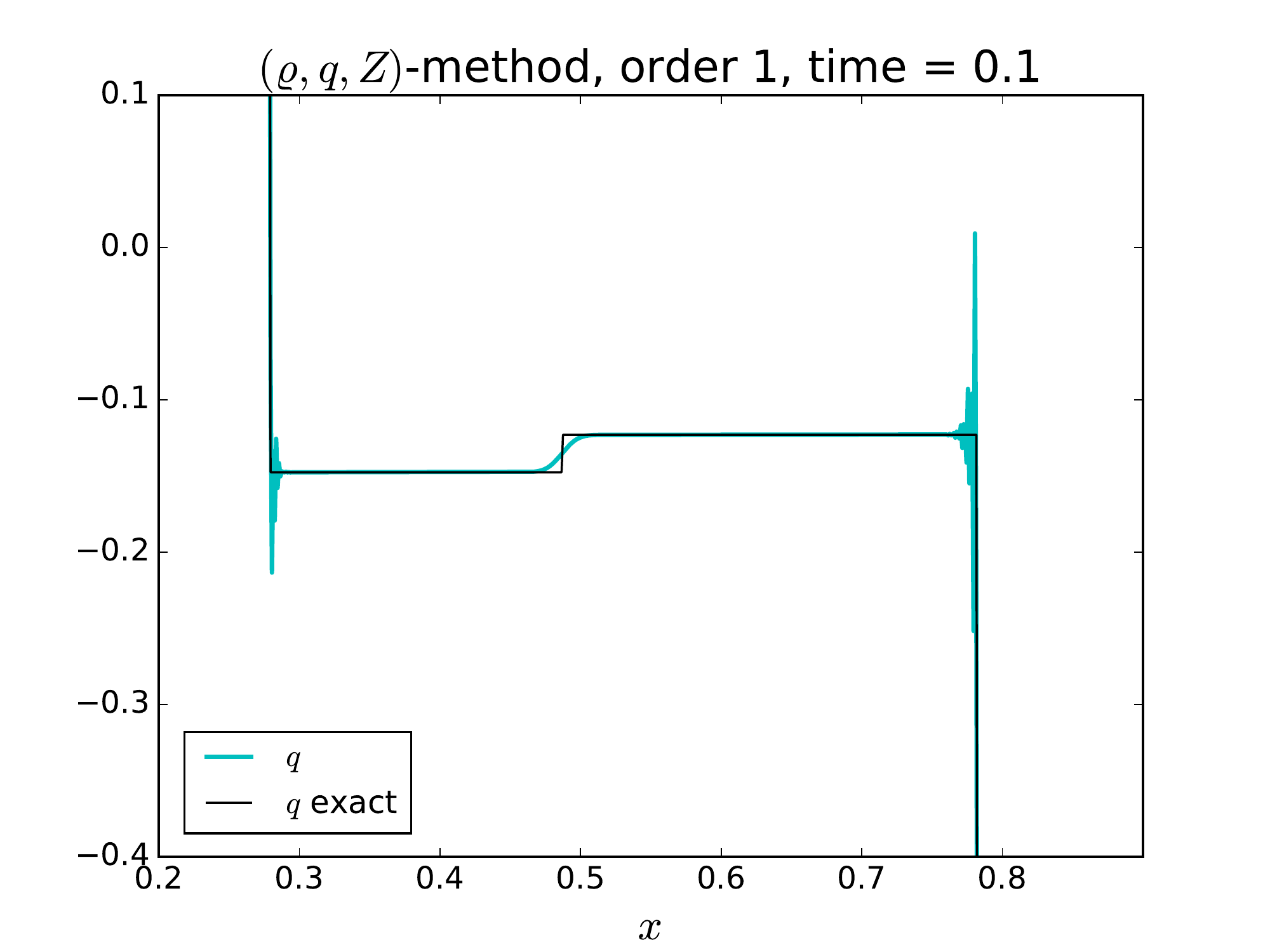}\includegraphics[width=0.45\textwidth]{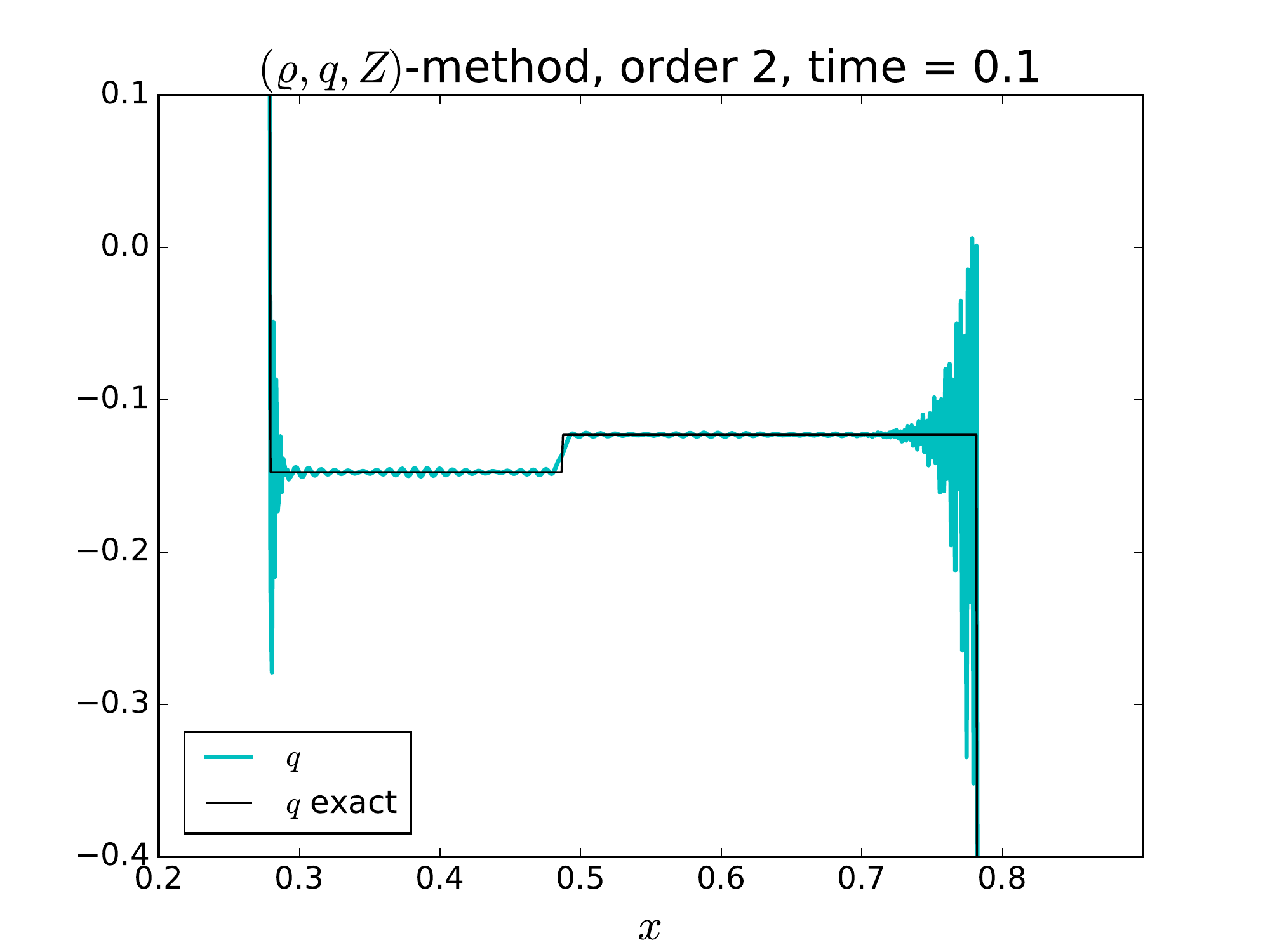}

\hspace{-0.6cm}\includegraphics[width=0.45\textwidth]{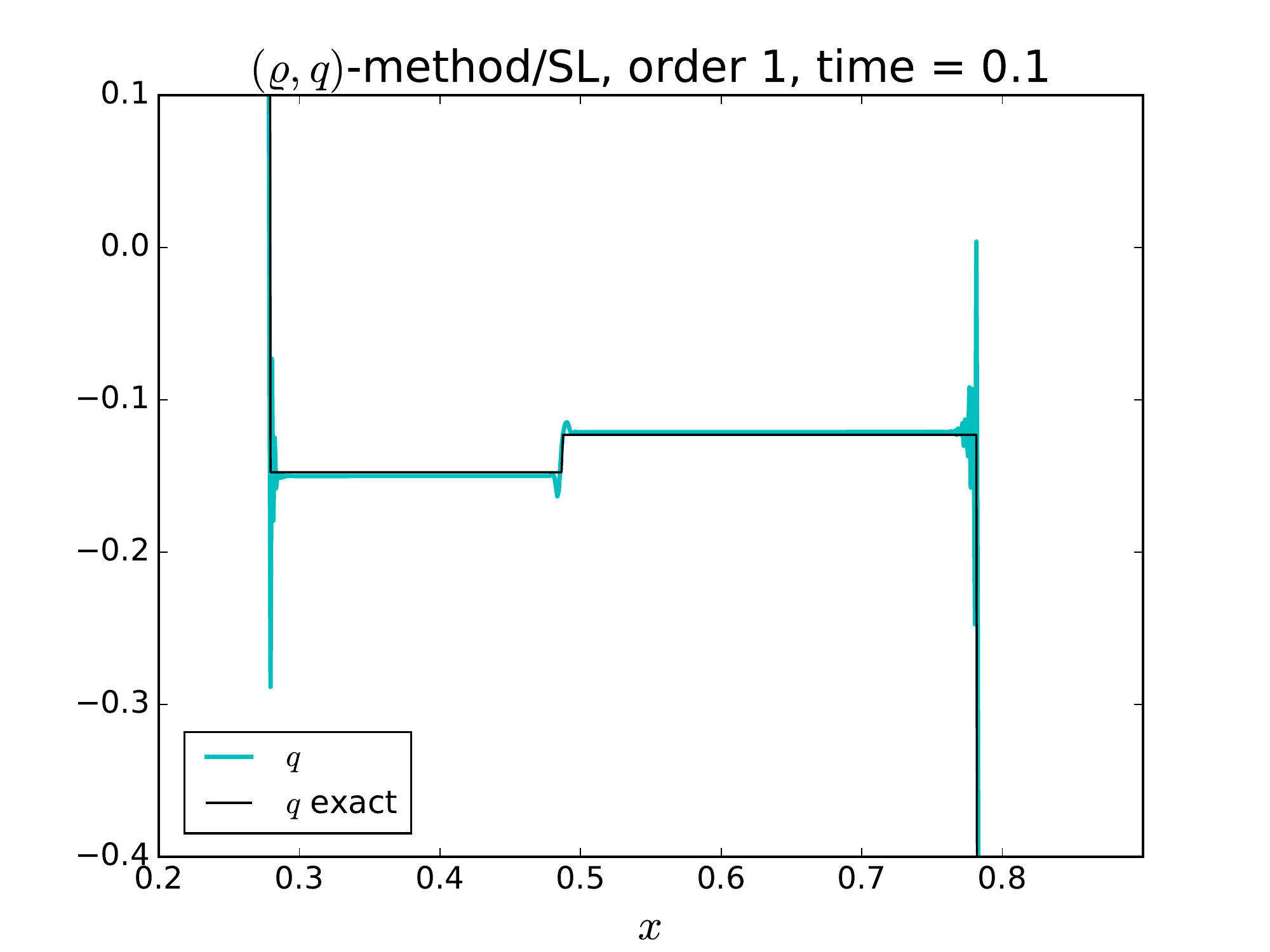}\includegraphics[width=0.45\textwidth]{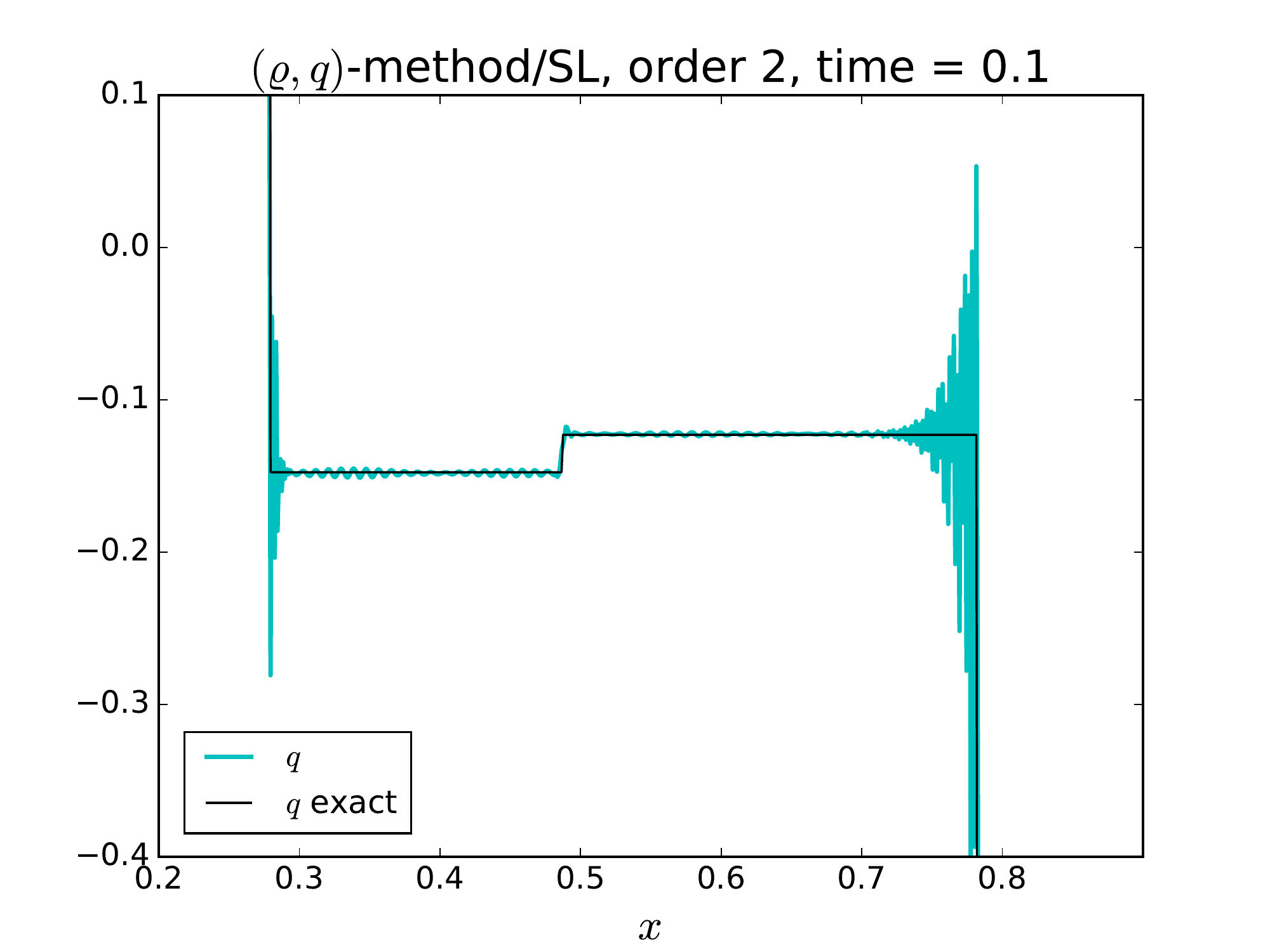}
\end{center}
\caption{Approximate momentum $q$ to Riemann problem \eqref{eq:riempb} at time $t = 0.1$ and comparison with the exact solution.  Numerical parameters: $\Delta x= 1\times 10^{-3}$, $ \Delta t = 0.1\, \Delta x$, $\alpha = 2$, $\gamma = 2$, $\eps=10^{-2}$.}
\label{sec:zoomriemtestcase}
\end{figure}

Oscillations in momentum are brought forth at the discontinuity interface of the shock waves. These oscillations are larger for second order schemes due to dispersion effects. In Figure \ref{sec:zoomriemtestcase}, we provide a zoom on these oscillations and compare the approximate solution to the exact one. The amplitudes of the oscillations are larger for the $(\vr, \q)$-method/SL  method. This may be the counterpart of the decoupling of the variables $(\vr,q)$ and $\vr^{\ast}$: in the computation of the implicit pressure (see eq.~\eqref{eq:ellipticRHOsemidiscrete}, left-hand side), $\vr$ and $\vr^{\ast}$ are not taken at the same time.

\paragraph{Stiff pressure ($\eps = 10^{-4}$)} With this value of $\eps$, the intermediate congested state has maximal wave speed equal to $\lambda_{\max} \approx 22$. Hence, taking time step $\Delta t$ equal to $0.1 \Delta x$ does not ensure the resolution of the fast waves. 

\begin{figure}
\begin{center}

\hspace{-0.6cm}\includegraphics[width=0.45\textwidth]{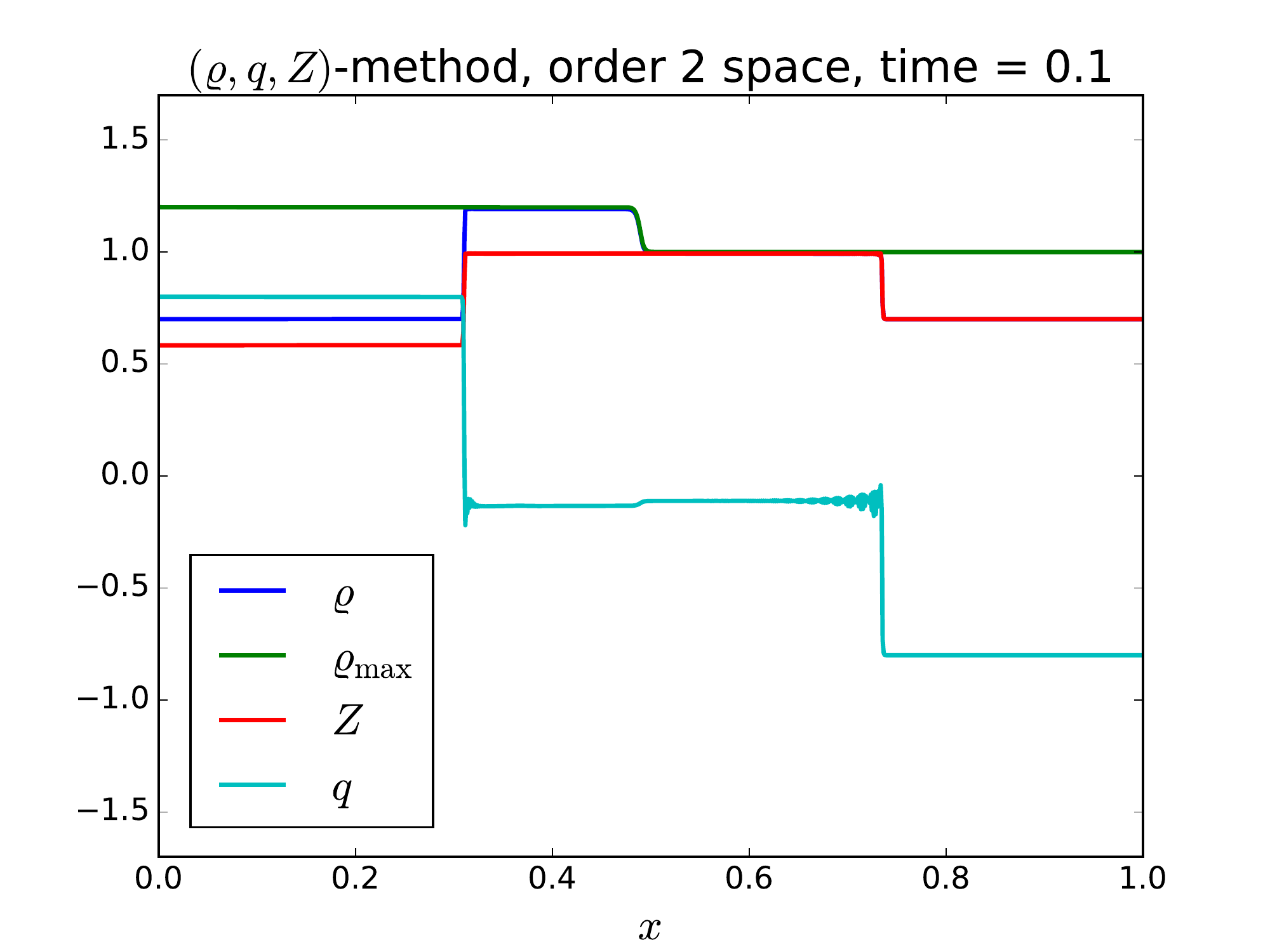}\includegraphics[width=0.45\textwidth]{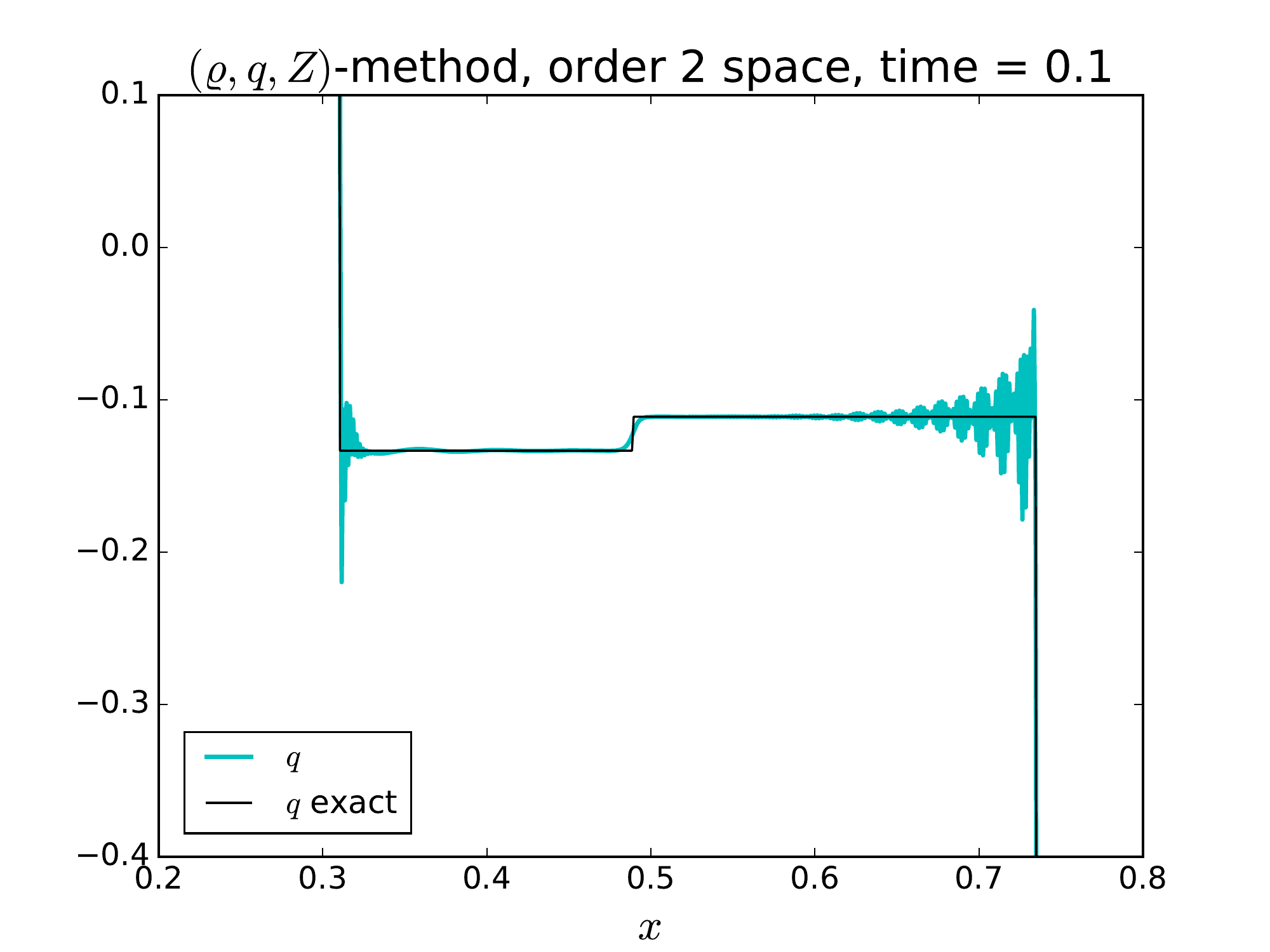}

\hspace{-0.6cm}\includegraphics[width=0.45\textwidth]{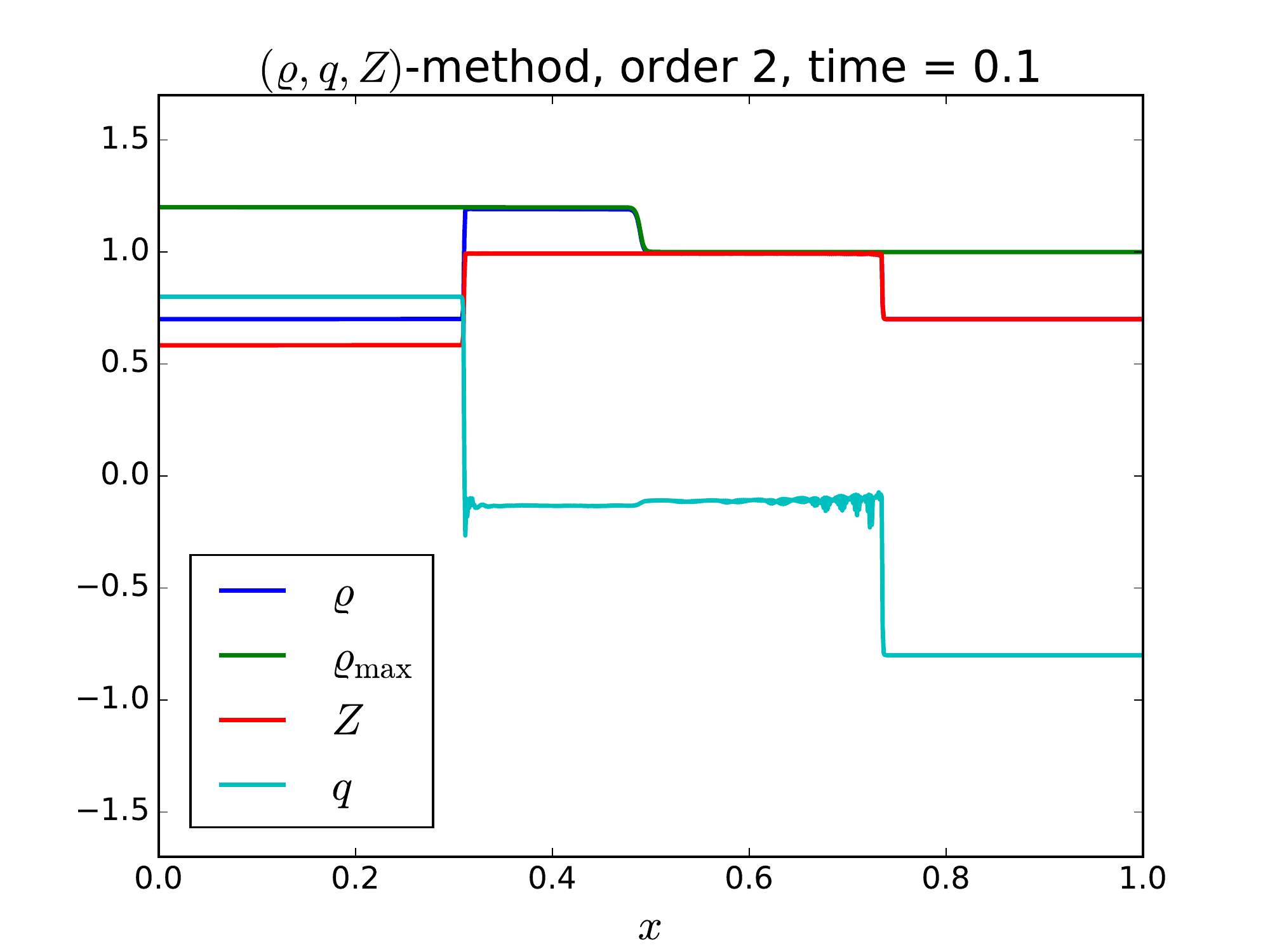}\includegraphics[width=0.45\textwidth]{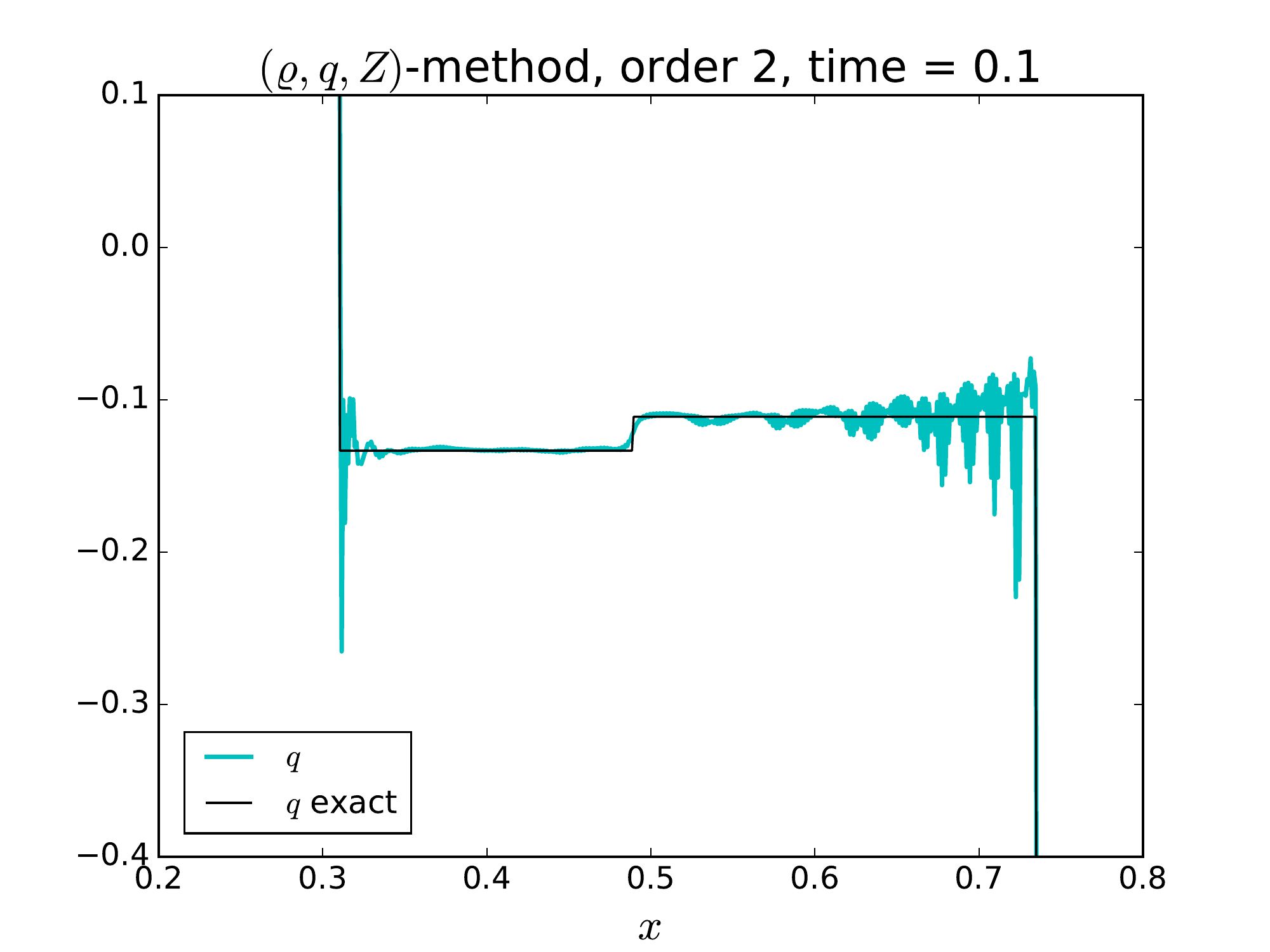}

\end{center}
\caption{Approximate solution to Riemann problem \eqref{eq:riempb} at time $t = 0.1$. Numerical parameters: $\Delta x= 1\times 10^{-3}$, $ \Delta t = 0.1\, \Delta x$, $\alpha = 2$, $\gamma = 2$, $\eps=10^{-4}$.}
\label{sec:riemtestcase_eps-4}
\end{figure}

Figure~\ref{sec:riemtestcase_eps-4} shows the solution at time $t=0.1$ using the $(\vr, \q, Z)$-method with second order in space accuracy. In the full second order scheme, the scheme switches automatically to a first order in time version of the scheme due to the large discontinuities in pressure, see Section \ref{sec:second_order_discretization}. We observe that the waves are well captured. As previously, oscillations in momentum develop at schock discontinuities and we observe that the second order in time version of the scheme leads to large uppershoots. In Table~\ref{tab:L1error}, we report the $L_{1}$ error between numerical and exact solution: we point out that the numerical errors are of the same order of magnitude independantly of the value of $\eps$. Quite similar results are obtained using the $(\vr, \q)$-method/SL. 

\begin{table}
\begin{tabular}{ll|c|c|c|c}
&&$\vr$&$q$&$Z$&$\vr^{\ast}$\\
\hline
$\eps = 10^{-2}$&order 2 in $x$&$8.66\times 10^{-4}$&$1.28\times 10^{-3}$&$3.03\times 10^{-4}$& $5.70\times 10^{-4}$\\
&order 2& $1.17\times 10^{-3}$ & $3.52\times 10^{-3}$&$5.89\times 10^{-4}$&$5.77\times 10^{-4}$\\
\hline
$\eps = 10^{-4}$&order 2 in $x$&$9.75\times 10^{-4}$&$2.11\times 10^{-3}$& $3.70\times 10^{-4}$& $5.71\times 10^{-4}$\\
&order 2& $9.89 \times 10^{-4}$& $3.04\times 10^{-3}$& $3.84\times 10^{-4}$& $5.77\times 10^{-4}$\\
\end{tabular}
\caption{$L_{1}$ error between the numerical solutions to Riemann problem \eqref{eq:riempb}  and exact solution at time $t = 0.1$. Numercial solution computed using the $(\vr,\q, Z)$-method. Numerical parameters: $\Delta x= 1\times 10^{-3}$, $ \Delta t = 0.1\, \Delta x$, $\alpha = 2$, $\gamma = 2$.}
\label{tab:L1error}
\end{table}

\subsection{Numerical convergence test-case}\label{ssec:num_conv}

\begin{figure}
\hspace{-0.6cm}\includegraphics[width=0.55\textwidth]{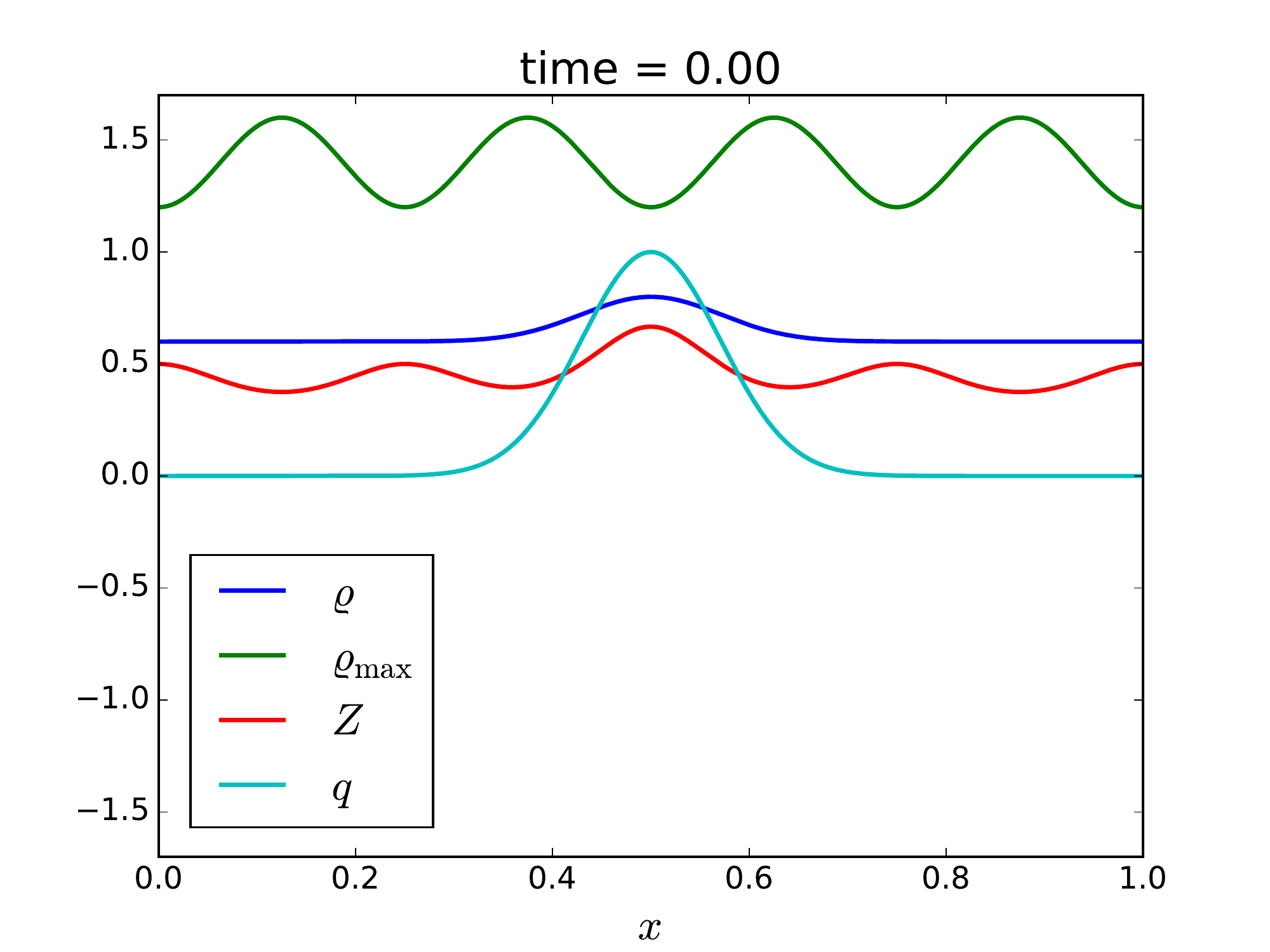}\includegraphics[width=0.55\textwidth]{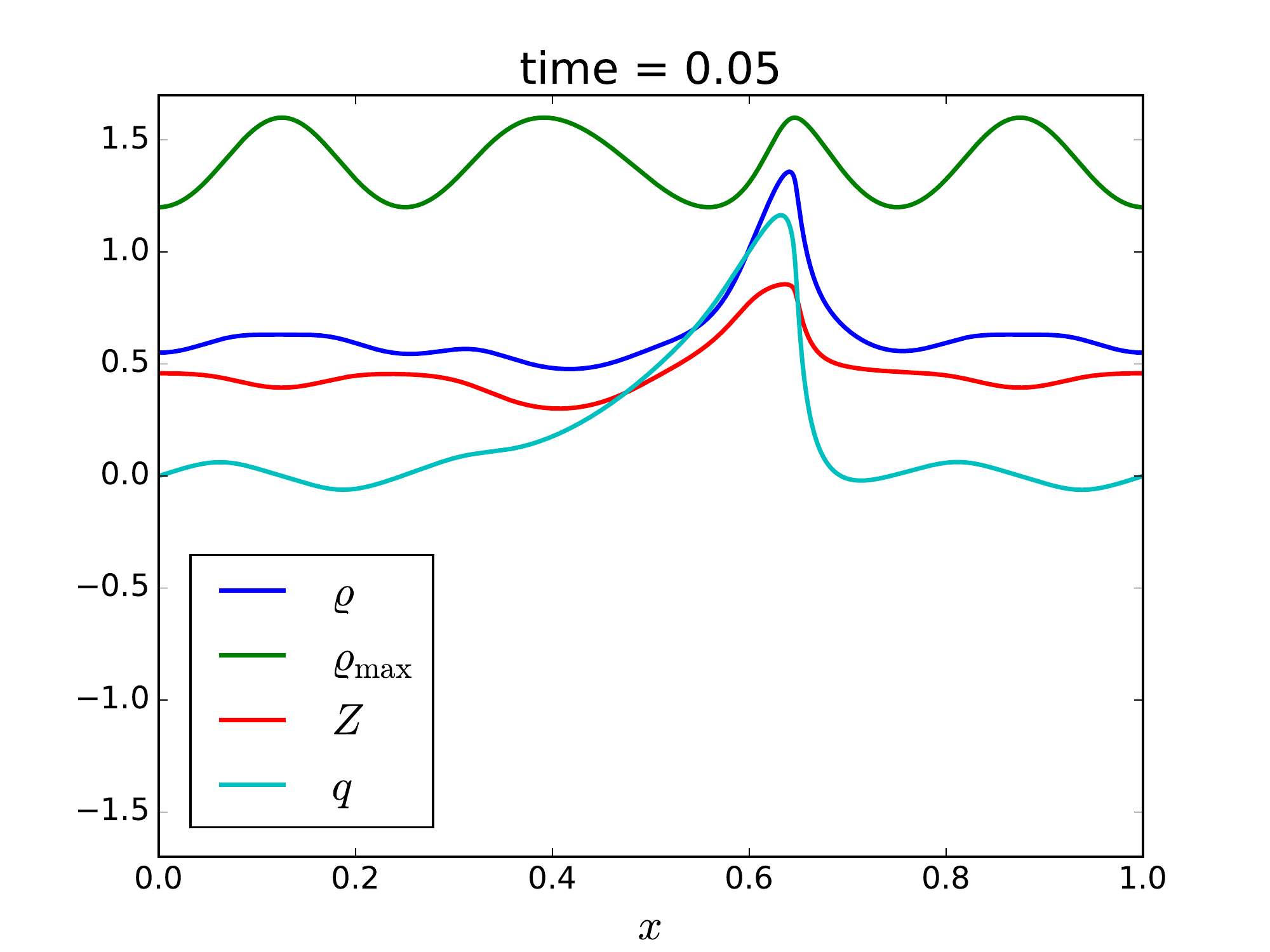}
\caption{Reference solution at initial time (left) and time $t=0.05$ (right). Numerical parameters: $\Delta x= 5\times 10^{-5}$, $ \Delta t = 0.1\, \Delta x$, $\gamma = 2$, $\eps=10^{-2}$.}
\label{sec:smoothtestcase}
\end{figure}

We here consider the following smooth initial data:
\begin{align*}
\vr_{0}(x) &= 0.6 + 0.2\, \exp\big(-(x-0.5)^2/0.01\big),\\
q_{0}(x) &= \exp\big(-(x-0.5)^{2}/0.01\big),\\
\vr^{\ast}_{0}(x) &= 1.2 + 0.2\, \big(1 - \cos\big(8 \pi (x-0.5)\big)\big),
\end{align*}
on the domain $[0,1]$ and perdiodic boundary conditions. We compute a reference solution at time $t = 0.05$ using the second order in space $(\vr, q, Z)$-method with small space and time steps $\Delta x= 5\times 10^{-5}$ and $ \Delta t = 0.1\, \Delta x$ (see Fig. \ref{sec:smoothtestcase}).

\begin{figure}
\hspace{-0.5cm}\includegraphics[width=1.2\textwidth]{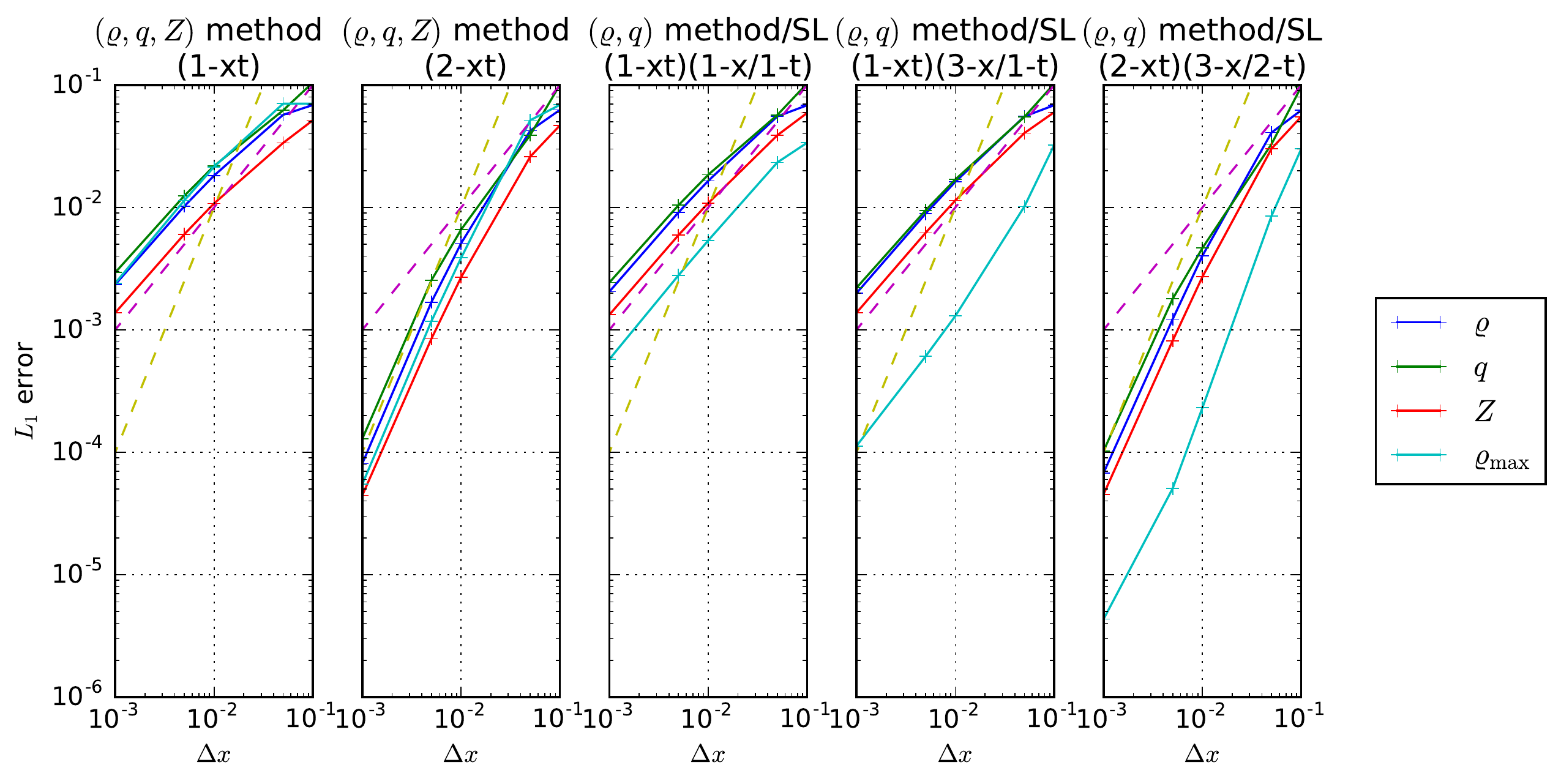}
\caption{$L_{1}$ errors for $\vr$, $q$, $Z$ and $\vr^{\ast}$ as function of $\Delta x$. Numerical parameters: $ \Delta t = 5\, \times 10^{-6}$ for first order scheme and $ \Delta t = 0.1\, \Delta x$ for second order scheme,  $\gamma = 2$, $\eps=10^{-2}$. $(\vr, \q, Z)$-method: ($k$-xt) $k$-th order in space and time. $(\vr,\q)$-method/SL: ($k$-xt)($m$-x/$n$-t) $k$-th order in space and time for the $(\vr,\q)$-method and $m$-th order in space and $n$-th order in time for the advection of $\vr^{\ast}$ by the semi-Lagrangian scheme. In dashed lines: first and second order curves.}
\label{sec:smooth-convergence}
\end{figure}

Figure \ref{sec:smooth-convergence} shows the $L_{1}$ errors between approximate solutions and the reference solution at time $t=0.05$ when the space step $\Delta x$ goes to $0$. For first order scheme, time step is set to $\Delta t = 5\times 10^{-6}$ while for second order schemes, time and space steps satisfy the relation $ \Delta t = 0.1\, \Delta x$ and both are varying. 

We observe that all the schemes exhibit their expected convergence rates. We point out that $(\vr, \q, Z)$-method and $(\vr, \q)$-method/SL have the same level of numerical errors except for variable $\vr^{\ast}$:  $\vr^{\ast}$ is better resolved with  $(\vr, \q)$-method/SL. This is all the more the case when using the third order semi-Lagrangian scheme (on the right two plots of Fig. \ref{sec:smooth-convergence}).

\section{Two-dimensional numerical results}\label{sec:2d_comp}
In this section we present the results of the numerical  simulations in two-dimensions. As for domain we take the unit square with the mesh size $\Delta x = 10^{-3}$ and the time-step $\Delta t = 10^{-4}$. In the following we choose singular pressure \eqref{pie} with the parameters $\ep = 10^{-4}$, $\alpha=2$, and the background pressure  \eqref{pbar} with the exponent $\gamma=2$, if not stated differently.

First part is devoted to comparison of $(\vr, \q, Z)$-method and $(\vr, \q)$-method/SL described in Section \ref{sec:numscheme}. Second is an application of $(\vr, q)$-method/SL to the evacuation scenario.  Third order in space semi-Lagrangian scheme is applied.

\subsection{Collision of 4 groups with variable congestion}\label{sec:test}

In the unit square periodic domain we specify 4 squares, with the centers in points $(x_c,y_c)=\{(0.2,0.5),(0.5,0.2),(0.5,0.8),(0.8,0.5) \}$. The length of the side $l$ of each square equals $0.2$ (for every square we introduce the notation $\text{Square}((x_c,y_c),l)$). We prescribe the initial momentum of $0.5$ pointing into the center of the domain provoking a collision. We consider three test cases varying in the initial congestion density, namely:
\begin{itemize}
\item[Case 1:] $\vr^*(x,0) = 1.0$;
\item[Case 2:] $\vr^*(x,0) = \begin{cases}
			 0.80  \text{ if } x\in \text{ Square}((0.2,0.5),0.2) \\
			 1.20  \text{ if } x\in \text{ Square}((0.5,0.2),0.2) \\
			 0.80  \text{ if } x\in \text{ Square}((0.8,0.5),0.2) \\
			 1.20  \text{ if } x\in \text{ Square}((0.5,0.8),0.2) \\
			 1.00 \text{ otherwise }
			\end{cases}
$;
\item[Case 3:] $\vr^*(x,0) = 1 + 0.05(\cos(10\pi x)+\cos(24\pi x))(\cos(6\pi y)+\cos(34 \pi y))$.
\end{itemize}

\begin{figure}
\includegraphics[width=1.\textwidth]{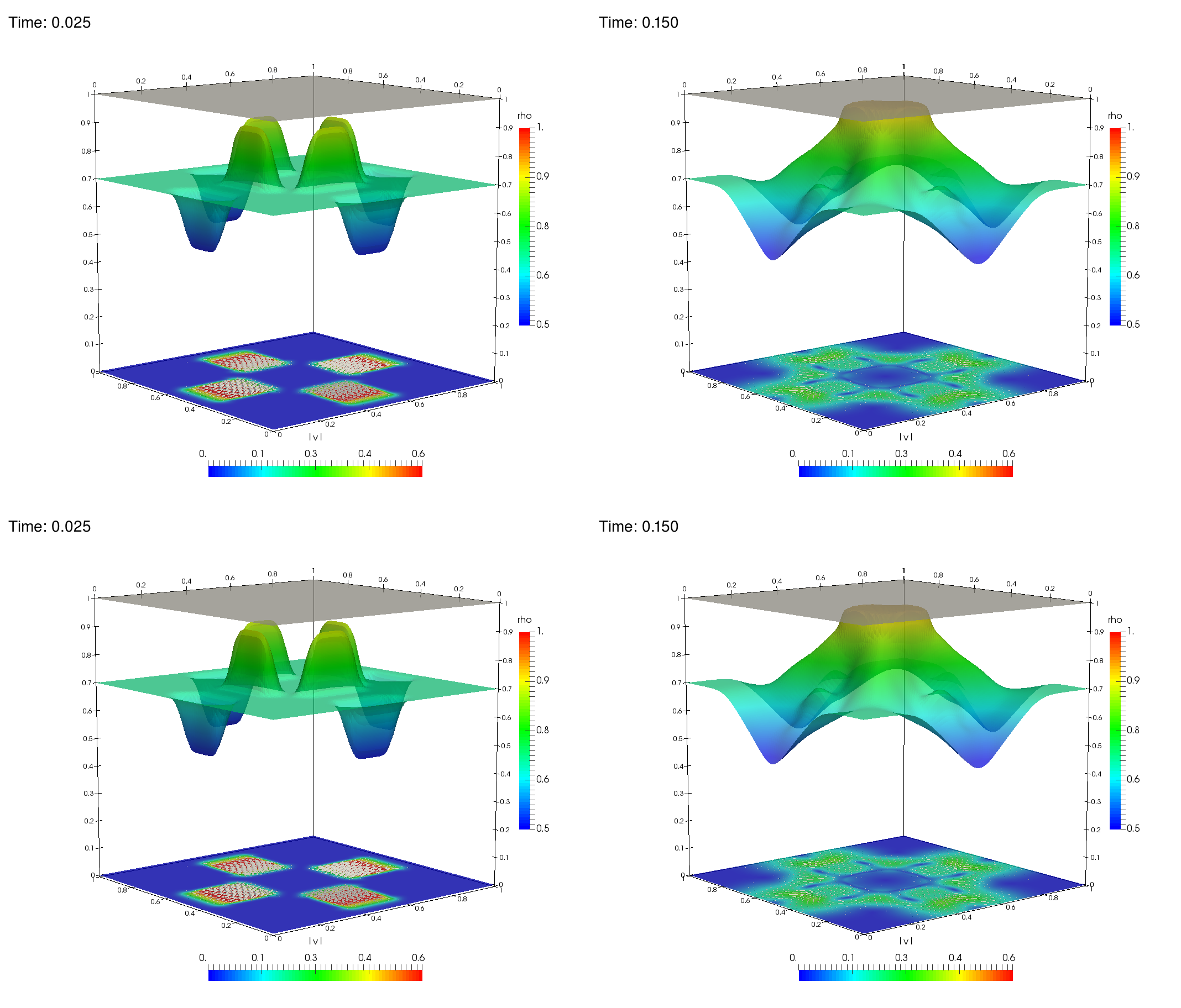}
\caption{ Case 1: the comparison of $(\vr, \q, Z)$-method (top) and $(\vr, \q)$-method/SL (bottom) at time $0.025$ (left), and $0.150$ (right).}
\label{fig:c1}
\end{figure}

\begin{figure}
\includegraphics[width=1.\textwidth]{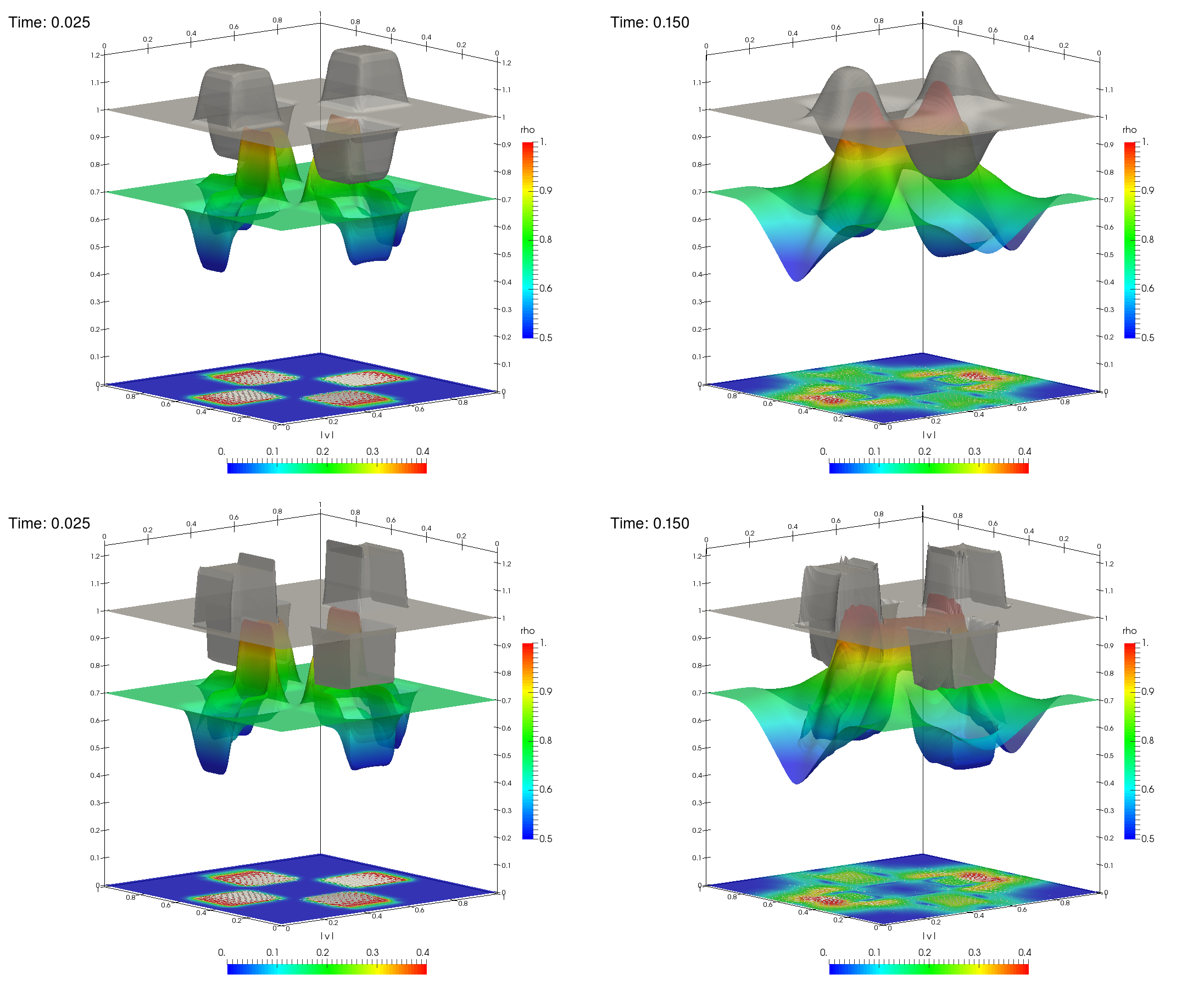}
\caption{ Case 2: the comparison of $(\vr, \q, Z)$-method (top) and $(\vr, \q)$-method/SL (bottom) at time $0.025$ (left), and $0.150$ (right).}
\label{fig:c2}
\end{figure}

\begin{figure}
\includegraphics[width=1.\textwidth]{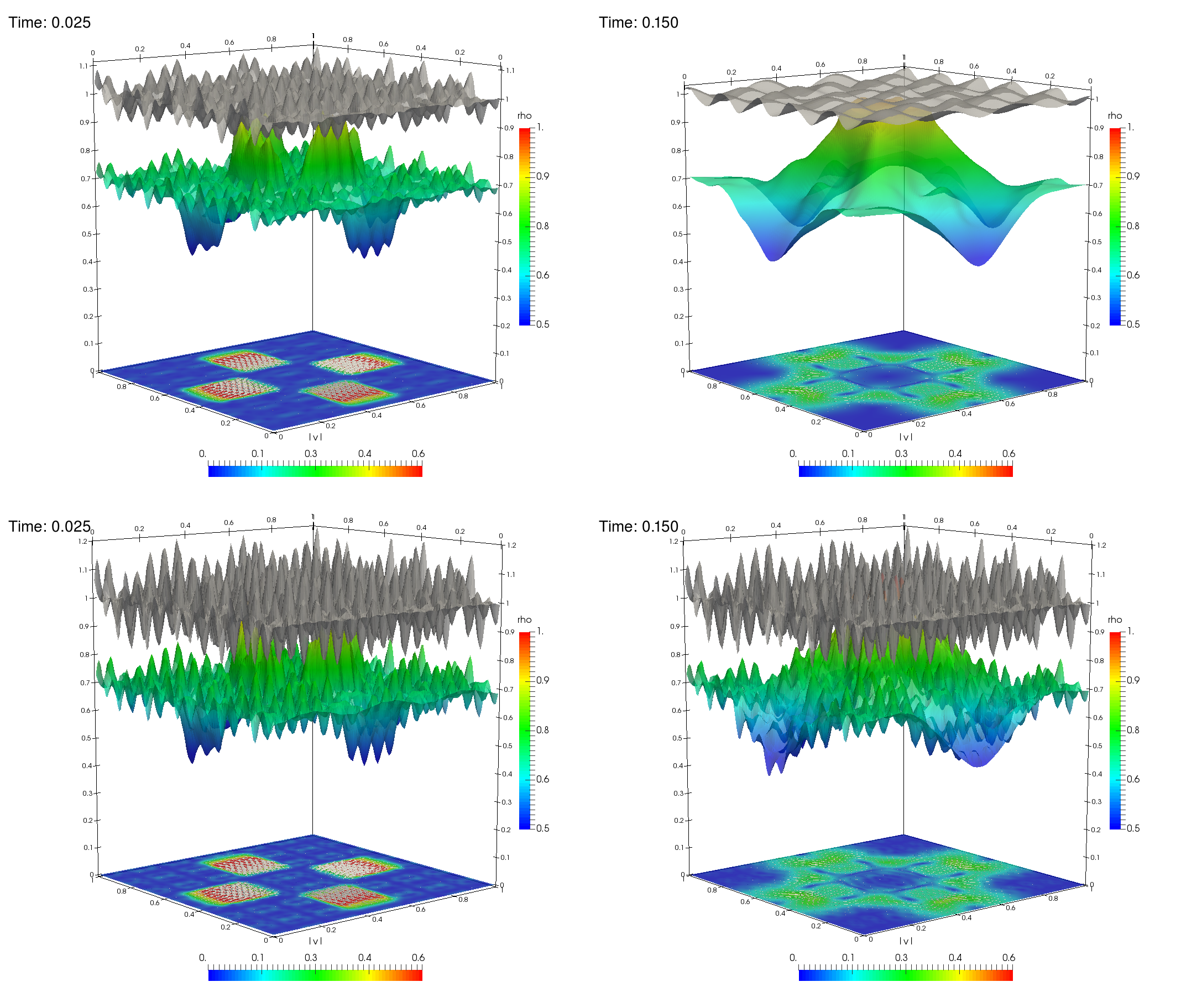}
\caption{ Case 3: the comparison of $(\vr, \q, Z)$-method (top) and $(\vr, \q)$-method/SL (bottom) at time $0.025$ (left), and $0.150$ (right).}
\label{fig:c3}
\end{figure}

The results of our simulations for these three cases are presented in Figures \ref{fig:c1}, \ref{fig:c2}, and \ref{fig:c3}, subsequently and in the Movies c1.mp4, c2.mp4, and c3.mp4. We see that in case of constant congestion density (Case 1, Figure \ref{fig:c1}, Movie c1.mp4) the two schemes provide almost identical outcome. The essential difference appears when $\vrs_0$ varies. We see in Figure \ref{fig:c2} (see also Movie c2.mp4) that the initial discontinuities of $\vrs$ are significantly smoothened by the $(\vr, \q, Z)$-method, while the $(\vr, \q)$-method/SL preserves the initial shape, which basically confirms our observations from Section \ref{ssec:num_conv}. This is even more visible in Figure \ref{fig:c3} (Movie c3.mp4), where the initial oscillations of $\vrs$ rapidly decay when simulated by the $(\vr, \q, Z)$-method.

Another interesting observation following from Figures \ref{fig:c2}, and \ref{fig:c3} (Movies c2.mp4, c3.mp4) when compared to Figure \ref{fig:c1} (Movie c1.mp4) is that the preference of the individuals $\vrs$ is significant factor to determine the density distribution even far away from the congestion zone. 

Moreover, comparing Figure  \ref{fig:c2} (Movie c2.mp4) with Figure  \ref{fig:c1} (Movie c1.mp4), we see a clear influence of the density constraint on the velocity of the agents. Indeed, for the Case 2, there is a significant disproportion between the velocities in the $x$ and $y$ directions at  time $t=0.150$  (see Figure \ref{fig:c2} right). This corresponds  to the fact that the agents moving toward the center along $y$ axis have 'more space' to fill  since  $\vrs$ for those groups is higher than the one for the groups moving in the $x$ direction. This results in a certain delay between collisions in two directions.

\subsection{Application to crowd dynamics}\label{sec:application}
In this section we investigate an influence of the variable density $\vr^*$ on a possible evacuation scenario. For this,  we consider an impenetrable room in the shape of unit square, initially filled with uniformly distributed agents. There is an exit located at $x \in [0.4,0.6],\  y = 0 $ that allows for free outflow. 
The initial density $\vr_0=0.6$ and the initial momentum is equal to $\vc{0}$.
The desire of going to the exit is introduced in the system \eqref{sysSl} \eqref{pie} by adding the relaxation therm in the momentum equation
  \eq{\label{relax}\partial_t \q + \Div \lr{\frac{\q\otimes \q}{\vr} + \pi_\ep\lr{\frac{\vr}{\vrs}}\Id + p\lr{\frac{\vr}{\vrs}}\Id}  = \frac{1}{\beta}\left(\q - \vr \boldsymbol{w}   \right),}
where $\boldsymbol{w}$ is the desired velocity, and $\beta$ stands for the relaxation parameter. The desired velocity is given by a unit vector field, that points into the centre of the exit, $\boldsymbol{w} = \left(-x/((x\!-\!0.5)^2\!+\!y^2) ,-y/((x-0.5)^2\!+\!y^2)\right)$ .

\begin{figure}
\begin{center}
\includegraphics[width=1.\textwidth]{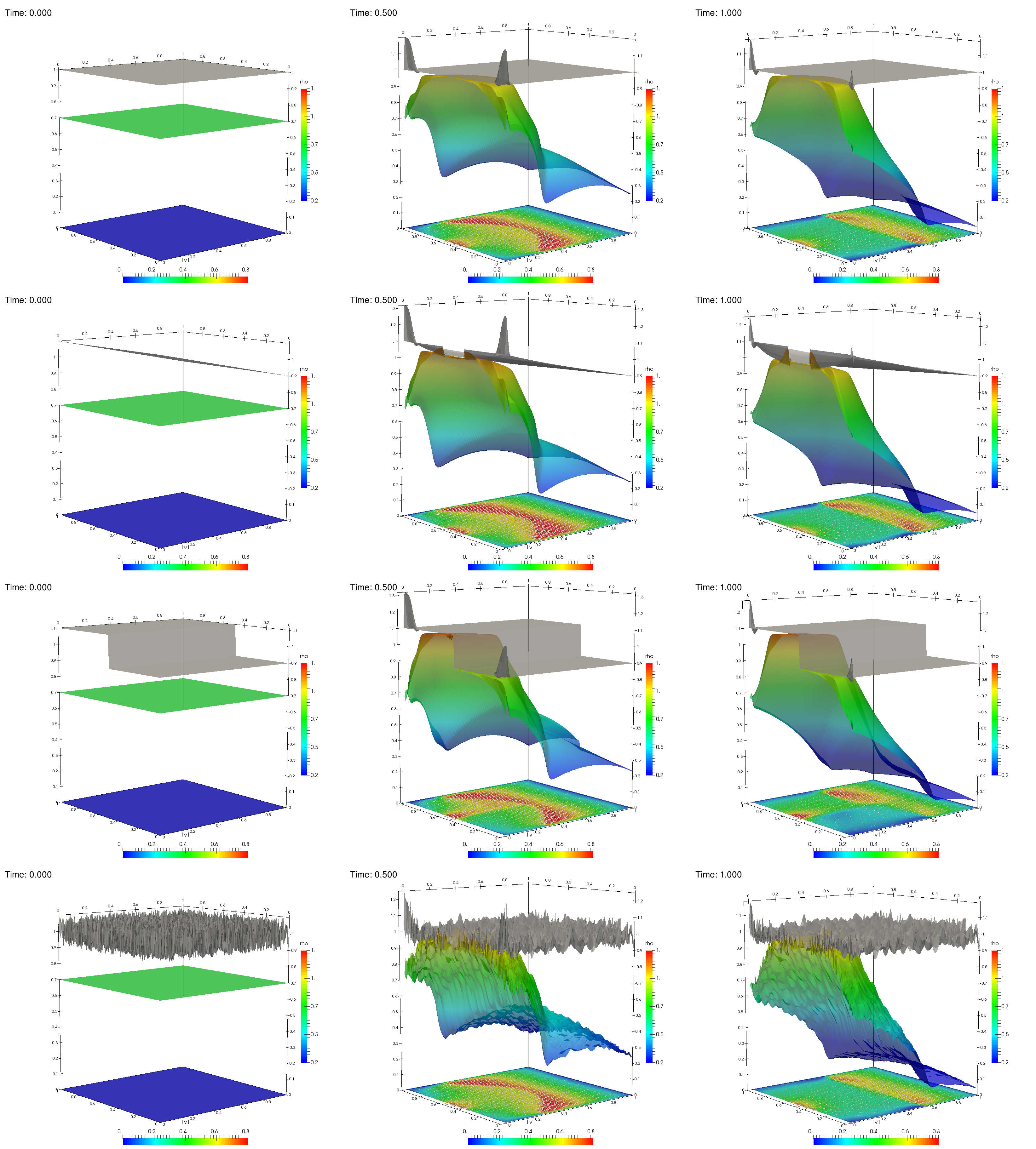}
\end{center}
\caption{Stop-and-go behaviour for the evacuation scenario, with $\vr^*_0$ being constant, with linear slope in $y-$direction \eqref{rhos02}, step-function \eqref{rhos0}, and a random function. The congestion density (upper) the density (middle) and the velocity amplitude (bottom) at times $t=0$ (left column) $t=0.5$ (middle column), and $t=1.0$ (right column).}
\label{fig:exit3d}
\end{figure}

\begin{figure}
\begin{center}
\includegraphics[width=1.\textwidth]{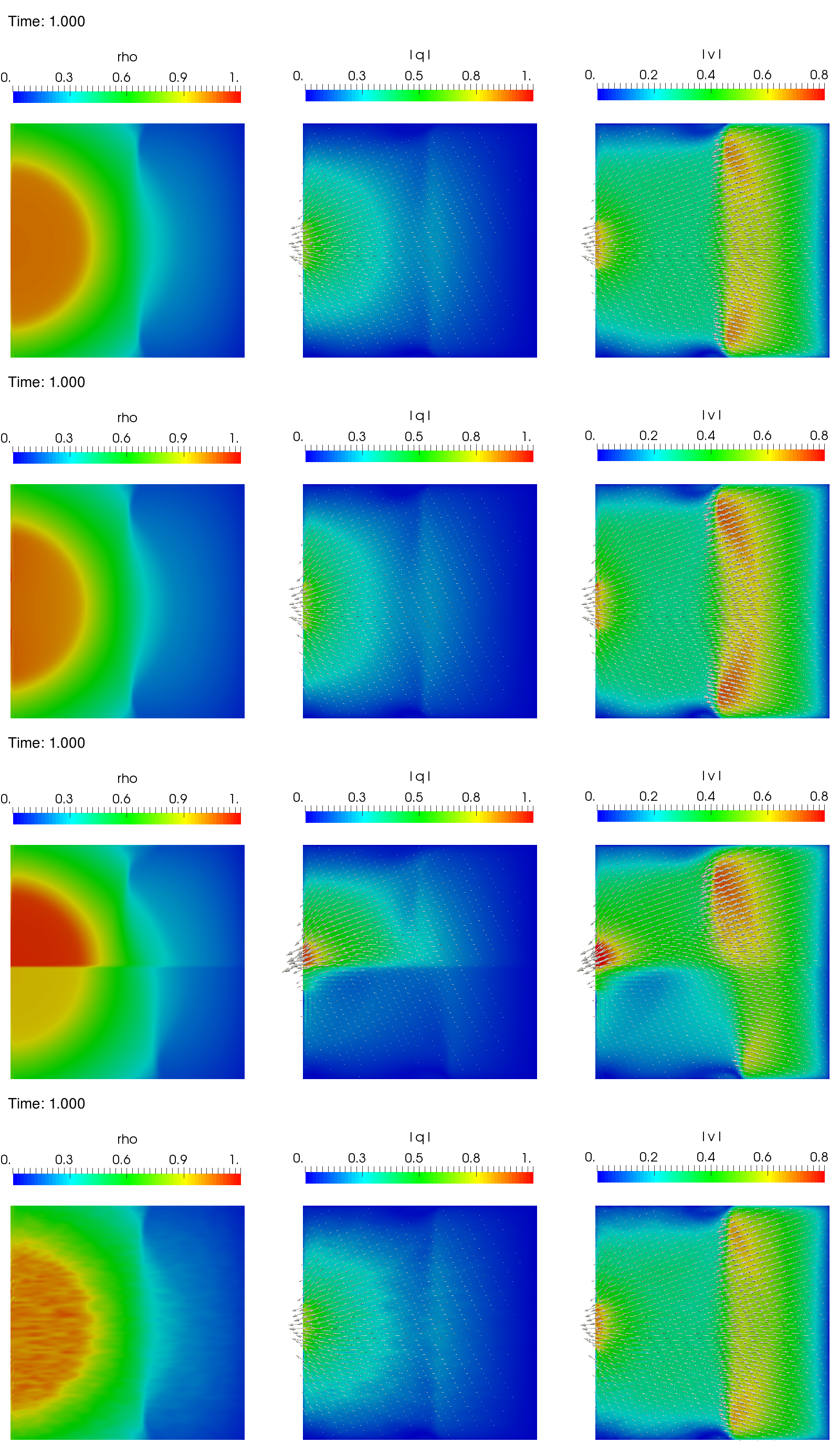}
\end{center}
\caption{The evacuation scenario for $\vr^*_0$ being constant, with linear slope in $y-$direction \eqref{rhos02}, step-function \eqref{rhos0}, and a random function. The figures present the values of the density $\vr$, the direction momentum $|q|$ and the direction and values of the velocity $v$ at time $t=1.0$ for different test cases.}
\label{exit_top}
\end{figure}

In the numerical scheme we apply splitting of the momentum equation between the transport and  pressure part, and the relaxation (source) part, with the intermediate momentum $\q^*$. After the momentum is updated we perform implicit relaxation step, for given density $\vr^{n+1}$,

\begin{subequations}
\begin{align}
&\frac{\q^* - \q^{n}}{\Delta t} + \nabla_{x}\cdot\left(\frac{\q^{n}\otimes \q^{n}}{\vr^{n}} + p\left(\frac{\vr^{n}}{{\vr^*}^n}\right) \Id \right) + \nabla_{x}\cdot\pi_\ep\left(\frac{\vr^{n+1}}{{\vr^*}^n}\right) = 0,\label{eq:semidiscreterelax1}\\
&\frac{\q^{n+1} - \q^*}{\Delta t} = \frac{1}{\beta}\left(\q^{n+1} - \vr^{n+1} \boldsymbol{w}   \right).\label{eq:semidiscreterelax2}
\end{align}
\label{eq:semidiscreterelax}
\end{subequations}
We use the $(\vr, \vc{q})$-method/SL, which requires to solve the transport equation for $\vrs$. This is especially problematic in the corners of the domain, where the Dirichlet boundary condition is considered. This leads to oscillations of $\vr$ and $\vrs$
close to these points.
Nevertheless, we may observe, see Figure \ref{fig:exit3d} and Figure \ref{exit_top} (see also Movies exit.mp4 and top.mp4), the so called {\it stop-and-go} behaviour, namely distinct high velocity regions in the domain, one in the vicinity of the exit and the second one that propagates in the direction opposite to flow. 
This reflects an empirical observation that  once a pedestrian arrives to the space of high congestion, he or she slows down or even stops until some space opens up in front. This kind of stop-and-go waves have been described, for example, by Helbing and Johansson in \cite{Helbing2011}. For the description of the real evacuation experiments we refer to \cite{HFV00}, see also \cite{garcimartin2014}. In the last of the mentioned papers the authors provide an experimental demonstration of the so called {\it faster goes slower} effect. This means that an increase in the density of pedestrians does not necessarily lead to a larger flow rate. 
Our simulations show that  when the parameter $\vr^*$ is low, the outflow of the individuals is slower. This is especially visible in the third row of  Figures \ref{fig:exit3d} and \ref{exit_top} presenting the evacuation scenario for the initial barrier density in the shape of the step function
\eq{\label{rhos0}
\vr^*_0(x,y)=\left\{\begin{array}{ccc}
1.1&\mbox{for} & 0.5<x<1,\\
0.9&\mbox{for} & 0<x< 0.5.
\end{array}
\right.}
This observation can be also confirmed in terms of speed of evacuation.

Indeed, we performed analogous simulations for 3 cases of constant $\vr^\ast_0$ equal to $0.9$, $1.0$. $1.1$ show that the speed of emptying the room is bigger the bigger value of $\vr^*_0$. To see this we have measured the mass remaining in the room at time $t=1$ and it is equal to $0.51030$, $0.048037$, and $0.457123$, respectively. 
We have moreover observed that evacuation speed of the room with individuals of the average  congestion preference equal to $1$ initially can be improved by placing the individuals with higher $\vrs_0$ closer to the exit. This is illustrated in the Figures \ref{fig:exit3d} and \ref{exit_top} the second row, for which, the initial congestion preference $\vr^\ast_0$ equals
\eq{\label{rhos02}
\vrs_0(x,y)=1.1-0.2y.
}
The random distribution of preferences of the individuals with expected value equal to $1$, on the other hand, corresponds to the increase of the evacuation time (see Figures \ref{fig:exit3d} and \ref{exit_top} the bottom row).

\section{Conclusion}

In this paper, we are interested in the numerical simulation of the Euler system with a singular pressure modeling variable congestion. As the stiffness of the pressure increases ($\varepsilon$ tends to $0$), the model tends to a free boundary  transition between compressible (non-congested) and incompressible (congested) dynamics. 

To numerically simulate the asymptotic dynamics, we propose an asymptotic preserving scheme based on a conservative formulation of the system and the methodology presented in \cite{DeHuNa}. We also propose a second order accuracy extension of the scheme following \cite{CorDeKum}. We then study the one-dimensional solutions to Riemann test-cases, their asymptotic limits and validate the code. We compare the results with those obtained with the scheme proposed in \cite{DeMiZa2016}. This latter scheme enables to better approximate the congestion density (at the contact wave) as soon as we use high accuracy in the advection of the congestion density. On the other hand, the former scheme seems to better preserve maximum principle on that variable. On two-dimensional simulations, we finally show the influence of this variable congestion density on the dynamics and show that  the model exhibit {\it stop-and-go} behavior.

The two schemes generate oscillations in momentum variable at discontinuities between congested and non-congested domain. This feature was already mentioned in \cite{DeHuNa}. This is all the more the case for the second order accuracy schemes. Specific method should be designed to cure this artefact.

\appendix
\section{Solution to the Riemann problem}
\label{sec:solRiemann}
The one-dimensional Riemann problem for the system \eqref{sysSZ} is the following initial-value problem:
\begin{subequations}
\begin{align}
&\partial_{t}\vr + \partial_{x}q = 0,\\
&\partial_{t}q + \partial_{x}\left(\frac{q^{2}}{\vr} + p_{\eps}(Z)\right) = 0,\\ 
&\partial_{t}Z + \partial_{x}\left(Z \frac{q}{\vr}\right) = 0,
\end{align}
\label{eq:system1D}
\end{subequations}
where $p_{\eps}(Z) = \pi_{\eps}(Z) + p(Z)$, and
\eq{\label{data1D}
(\vr,q,Z)(0,x)=\left\{
\begin{array}{lll}
(\vr_\ell,q_\ell,Z_\ell)&\mbox{for}&x<0,\\
(\vr_r,q_r,Z_r)&\mbox{for}&x>0.
\end{array}\right.
}
The purpose of this section is to find possible weak solution to \eqref{eq:system1D} \eqref{data1D}. We will also consider the limit of these solutions   as $\eps \rightarrow 0$. 

As already mentioned in the introduction, the system \eqref{eq:system1D} is strictly hyperbolic provided $p_\ep'(Z)>0$, see \eqref{eq:eigenval}. The associated characteristic fields are given by:
\begin{equation*}
r_{1}^{\eps}(\vr,q,Z) = \begin{bmatrix}
1\\
v-\displaystyle\sqrt{\frac{Z}{\vr} p_{\eps}'(Z)}\\
Z/\vr
\end{bmatrix},\ 
r_{2}^{\eps}(\vr,q,Z)  = \begin{bmatrix}
1\\
v\\
0
\end{bmatrix},\ r_{3}^{\eps}(\vr,q,Z)  = \begin{bmatrix}
1\\
v+\displaystyle\sqrt{\frac{Z}{\vr} p_{\eps}'(Z)}\\
Z/\vr
\end{bmatrix},
\end{equation*}
where $v = q/\vr$ is the velocity. The second characteristic field is linearly degenerate (since $\nabla \lambda_{2}\cdot r_{2} = 0$).  The two others characteristic field are genuinely non-linear. 

We now present the elementary wave solutions of the Riemann problem.

\subsection{Elementary waves}

\paragraph{Shock discontinuities} A shock wave is a discontinuity between two constant states, $(\vr, q, Z)$ and $(\hat\vr, \hat q, \hat Z)$, travelling at a constant speed $\sigma$. We now fix the left (or right) state $(\hat\vr, \hat q, \hat Z)$ and look for all triples $(\vr, q, Z)$ that can be connected to $(\hat\vr, \hat q, \hat Z)$ by the shock discontinuity. Across the shock, the  Rankine-Hugoniot conditions must be satisfied meaning that:
\begin{equation*}
[q] = \sigma [\vr],\quad \left[\frac{q^{2}}{\vr} + p_{\eps}(Z)\right] = \sigma [q],\quad \left[Z \frac{q}{\vr}\right] = \sigma [Z],
\end{equation*}
where $[a] := a - \hat a$ denotes the jump of quantity $a$. Treating $\vr$ as a parameter, we check that the two admissible states are of the form $(\vr, q_{h,\pm}(\vr), Z(\vr))$ with $q_{h,\pm}=\vr v_{h,\pm}(\vr)$ and
\eqh{
&v_{h,\pm}(\vr) = \hat v \pm \text{sign}(Z(\vr)-\hat Z)\frac{1}{\sqrt{\hat\vr\vr}}\sqrt{(\vr - \hat\vr)\lr{p_{\eps}\lr{\frac{\hat Z \vr}{\hat \vr}} - p_{\eps}(\hat Z)}},\\
&Z(\vr) = \hat Z\frac{\vr}{\hat \vr}.
}
The shock speed therefore equals:
\begin{equation*}
\sigma_{\pm} = \hat v \pm  \text{sign}(Z-\hat Z)\sqrt{\frac{\vr}{\hat\vr}}\sqrt{\frac{p_{\eps}(\hat Z \vr/\hat \vr) - p_{\eps}(\hat Z)}{(\vr - \hat\vr)}}.
\end{equation*}
These solutions can also be expressed as  functions of $Z$:
\eqh{
&\vr(Z) = Z\frac{\hat \vr}{\hat Z},\\
&v_{h,\pm}(Z) = \hat v \pm  \text{sign}(Z-\hat Z) \frac{1}{\sqrt{\hat\vr}}\sqrt{\lr{1 - \frac{\hat Z}{Z}} \big(p_{\eps}(Z) - p_{\eps}(\hat Z)\big)}.
}
Note that the maximal density ($\vr^{\ast} = \vr/Z$) does not jump across a shock discontinuity. 
Expanding $(\vr(Z), q_{h,\pm}(Z), Z)$ around $Z = \hat Z$, we obtain
\begin{align*}
\vr(Z) - \hat\vr &= (Z-\hat Z)\frac{\hat \vr}{\hat Z},\\
\vr(Z)v_{h,\pm}(Z) - \hat \vr\hat v &= (Z-\hat Z)\frac{\hat \vr}{\hat Z}\hat v  \pm  Z\frac{\hat \vr}{\hat Z}\text{sign}(Z-\hat Z) \sqrt{\frac{1}{\hat\vr}}\sqrt{(1 - \hat Z/Z) \big(p_{\eps}(Z) - p_{\eps}(\hat Z)\big)}\\
&\approx (Z-\hat Z)\frac{\hat \vr}{\hat Z}\hat v  \pm  Z\frac{\hat \vr}{\hat Z}\text{sign}(Z-\hat Z) \sqrt{\frac{1}{\hat\vr Z}}\sqrt{p_{\eps}'(\hat Z) (Z-\hat Z)^{2}}\\
&\approx (Z-\hat Z)\frac{\hat \vr}{\hat Z} \left(\hat v  \pm  \sqrt{\frac{\hat Z}{\hat\vr}}\sqrt{p_{\eps}'(\hat Z)}\right),\\
Z - \hat Z &= (Z-\hat Z)\frac{\hat \vr}{\hat Z} \frac{\hat Z}{\hat \vr}. 
\end{align*}
Note that $(\vr(Z),q_{h,-}(Z), Z)$ is tangent at $(\hat\vr,\hat q,\hat Z)$ to $r_1(\hat\vr,\hat q,\hat Z)$, therefore $v_{h,-}$ corresponds to the $1$-characteristic field, analogously $v_{h,+}$ corresponds to the $3$-characteristic field. The graph of $Z \mapsto v_{h,-}(Z)$ (resp. $Z \mapsto v_{h,+}(Z)$) is called the 1-Hugoniot curve (resp. 3-Hugoniot curve) issued from $(\hat v, \hat Z)$.

To check the admissibility of the discontinuity, we need to check the entropy condition. If $(\hat v, \hat Z)$ is the left state, the right states that can be connected to it by an entropic shock wave are those located on the 1-shock curve $\left\{ \lr{v_{h,-}(Z), Z}: Z > \hat Z\right\}$ or the 3-shock curve $\left\{ (v_{h,+}(Z), Z):Z < \hat Z\right\}$. Indeed, on these curves the associated eigenvalue is decreasing. If on the other hand, $(\hat v, \hat Z)$ is the right state, the left states that can be connected to it by an entropic shock wave are those located on the 1-shock curve $\left\{ (v_{h,-}(Z), Z): Z < \hat Z\right\}$ or the 3-shock curve $\left\{ (v_{h,+}(Z), Z): Z > \hat Z\right\}$. Indeed, on these curves the associated eigenvalue is increasing.

\paragraph{Rarefaction waves} The rarefaction waves are continuous self-similar solutions, $(\vr(t,x), q(t,x), Z(t,x))=(\vr(x/t), q(x/t), Z(x/t))$, connecting two constant states $(\vr, q, Z)$ and $(\hat\vr, \hat q, \hat Z)$. They thus satisfy the following differential equations:
\begin{equation}\label{ssym}
\vr'(s) = 1,\quad q'(s) = \tilde v(s)\pm\displaystyle\sqrt{\frac{Z(s)}{\vr(s)} p_{\eps}'(Z(s))},\quad Z'(s) = Z(s)/\vr(s),
\end{equation}
Denoting $q(s) = \vr(s) \tilde v_{i,\pm}(s)$ and parametrizing by $\vr$, we obtain:
\begin{equation*}
\tilde v_{i,\pm}'(\vr) = \pm \frac{1}{\vr}\displaystyle\sqrt{\frac{Z(\vr)}{\vr} p_{\eps}'(Z(\vr))},\quad Z'(\vr) = Z(\vr)/\vr.
\end{equation*}
From the first and third equation of \eqref{ssym}, we have $(\vr/Z(\vr))' = 0$, and so,  $\vr/Z(\vr) = \hat\vr/Z(\hat \vr)$.
This means that as in the case of shock discontinuities the maximal density $\vrs$ does not jump.  Denoting $\vr^{*} = \vr/Z(\vr)$ and making the change of coordinates $v_{i,\pm}(Z) = \tilde v_{i,\pm}(\vr)$ with $\vr = \vr^{\ast} Z$, we thus have:
\begin{equation*}
v_{i,\pm}'(Z) = \pm \frac{1}{Z}\displaystyle\sqrt{\frac{1}{\vr^{\ast}} p_{\eps}'(Z)}.
\end{equation*}
Hence, the states satisfy:
\begin{equation}
v_{i,\pm} (Z) = \hat v \pm  \big(F_{\eps}(Z) - F_{\eps}(\hat Z)\big),
\end{equation}
where $F_{\eps}$ is an antiderivative of $Z \mapsto \frac{1}{ Z}\sqrt{\frac{1}{\vrs} p_{\eps}'(Z)}$. 

 The graph of $Z \mapsto v_{i,+}(Z)$ (resp. $Z \mapsto v_{i,-}(Z)$) is called the 1-integral curve (resp. 3-integral curve) issued from $(\hat v, \hat Z)$. If $(\hat v, \hat Z)$ is a left state, the right states that can be connected to it by an entropic rarefaction wave are those located on the 1-integral curve $\left\{ (v_{i,-}(Z), Z):Z < \hat Z\right\}$ or the 3-integral curve $\left\{ (v_{i,-}(Z), Z): Z > \hat Z\right\}$. Indeed, on these curves the associated eigenvalue is increasing. If $(\hat v, \hat Z)$ is a right state, the left states that can be connected to it by an entropic rarefaction wave are those located on the 1-integral curve $\left\{ (v_{i,-}(Z), Z):Z > \hat Z\right\}$ or the 3-integral curve $\left\{ (v_{i,-}(Z), Z): Z < \hat Z\right\}$. Indeed, on these curves the associated eigenvalue is decreasing.

\paragraph{Contact discontinuities} Since the second characteristic field is linearly degenerate, there are linear discontinuities that propagate at velocity $\lambda_{2} = \hat v$. Let us write the Rankine-Hugoniot conditions:
\begin{equation*}
[q] = \hat v [\vr],\quad \left[\frac{q^{2}}{\vr} + p_{\eps}(Z)\right] = \hat v [q],\quad \left[Z \frac{q}{\vr}\right] = \hat v [Z].
\end{equation*}
From the first relation, we obtain $v = \hat v$ and then the second relation states that the pressure jump is zero. By strict monotony of the pressure, it implies that $Z = \hat Z$ and the third equation is satisfied. Along this discontinuity, the velocity and the pressure are thus conserved. Note that every density jump is possible. 

\subsection{Solution to Riemann problem}

Let $(\vr_{\ell}, q_{\ell}, Z_{\ell})$ and $(\vr_{r}, q_{r}, Z_{r})$ be the left and right initial states \eqref{data1D}. The solutions to Riemann problems are determined as follows. First,  in the $(v, Z)$ plane, find out the intersection state $(v_{m}, Z_{m})$ of the 1-st integral/Hugoniot curves issued from $(v_{\ell}, Z_{\ell})$ and the 3-rd integral/Hugoniot curves issued from $(v_{r}, Z_{r})$. Then, compute the two densities $\vr_{m, \ell}$ and $\vr_{m, r}$ so that the congestion density across the two non-linear waves is conserved. Then we connect the two distinct intermediate states by a contact discontinuity. We finally end up with the following solution: 
\begin{align}
(\vr_{\ell}, q_{\ell}, Z_{\ell})\ \overset{shock/rarefaction}{\rightarrow} &(\vr_{m, \ell}, \vr_{m, \ell}v_{m}, Z_{m})\nonumber\\
\overset{contact}{\rightarrow}\ &(\vr_{m, r}, \vr_{m, r}v_{m}, Z_{m})\   \overset{shock /rarefaction}{\rightarrow} (\vr_{r}, q_{r}, Z_{r}) 
\label{eq:solutionRiemPb}
\end{align}
where $\vr_{m, \ell} = Z_{m} \vr_{\ell}/Z_{\ell}$ and $\vr_{m, r} = Z_{m} \vr_{r}/Z_{r}$. The nature of the non-linear waves (rarefaction or shock) depends on the relative position of the states $(v_{\ell}, Z_{\ell})$, $(v_{r}, Z_{r})$ in the $(v, Z)$ plane.

\subsection{Limit $\eps \rightarrow 0$}
We are now interested in the asymptotic behaviour, when $\eps\to0$ of the Hugoniot $v_{h,\pm}^\eps$ and the integral curves $v_{i,\pm}^{\eps}$ obtained in the previous paragraph for the elementary waves. We have the following result.
\begin{prop} The graph of the Hugoniot curve, $\left\{ (Z, v_{h,\pm}^{\eps}(Z)) :Z \in [0,1)\right\}$, tends to the union of the set $\left\{ (Z, v_{h,\pm}^{0}(Z)): Z \in [0,1)\right\}$ and the horizontal half straight line $\left\{ (1, v):  v \in [v_{h,\pm}^{0}(1), +\infty)\right\}$. 

The graph of the integral curve, $\left\{ (Z, v_{i,\pm}^{\eps}(Z)):  Z \in [0,1)\right\}$, tends to the union of the set $\left\{ (Z, v_{i,\pm}^{0}(Z)):  Z \in [0,1)\right\}$ and the horizontal half straight  line $\left\{ (1, v): v \in [v_{i,\pm}^{0}(1), +\infty)\right\}$. 
\end{prop}
\noindent The proof of this proposition uses the convexity of the pressure and are similar to the one developed in \cite{DeNaBonSan}.

Regarding the Riemann problem in the limit $\eps \rightarrow 0$, the intersection  point of the 1-st integral/Hugoniot curves issued from $(v_{\ell}, Z_{\ell})$ and the 3-rd integral/Hugoniot curves issued from $(v_{r}, Z_{r})$, denoted by $(v_{m}^{\eps}, Z_{m}^{\eps})$, has either a limit $(v_{m}^{0}, Z_{m}^{0})$ with $0 \leqslant Z_{m}^{0} < 1$ or tends to a congested state $(\bar v,1)$. Then finding a solution can be divided into the following steps:
\begin{itemize}
\item[(1)] compute the intersection $(v_{m}^{0}, Z_{m}^{0})$ of the 1-st integral/Hugoniot curves and 3-rd integral/Hugoniot curves;
\item[(2a)] if $Z_{m}^{0} < 1$, the solution is as described in the previous section, it is a usual Riemann solution of the hyperbolic system with no congestion pressure;
\item[(2b)] if $Z_{m}^{0} \geq1$, then the congested state is given by the following proposition.
\end{itemize}
\begin{prop}[Case $Z_{m}^{0} \geq 1$.]  The solution consists in three waves:
\begin{equation*}
(\vr_{\ell}, q_{\ell}, Z_{\ell})\ \overset{shock}{\rightarrow}\ (\vr_{\ell}^{\ast}, \vr_{\ell}^{\ast}\bar v, 1)\  \overset{contact}{\rightarrow}\ (\vr^{\ast}_{r}, \vr^{\ast}_{r}\bar v, 1)\   \overset{shock}{\rightarrow} (\vr_{r}, q_{r}, Z_{r}) 
\end{equation*}
where the intermediate velocity $\bar v$ and pressure $\bar p$ satisfy:
\begin{align*}
&\bar v  =  v_{\ell} -   \sqrt{\frac{1}{\vr_{\ell}}}\sqrt{(1 - Z_{\ell}) (\bar p - p_{0}(Z_{\ell}))} = v_{r} +   \sqrt{\frac{1}{\vr_{r}}}\sqrt{(1 - Z_{r}) (\bar p - p_{0}(Z_{r})) },
\end{align*}
the intermediate densities are given by:
\begin{equation*}
\hat\vr_{\ell} = \vr_{\ell}/Z_{\ell} = \vr^{\ast}_{\ell},\quad \hat\vr_{r} = \vr_{r}/Z_{r} = \vr^{\ast}_{r},
\end{equation*}
and the shock speeds $\sigma_{-}$, $\sigma_{+}$ are given by:
\begin{equation*}
\sigma_{-} = v_{\ell} - \sqrt{\frac{\vrs_{\ell}}{\vr_{\ell} (\vrs_{\ell} - \vr_{\ell})}}\sqrt{\bar p - p_{0}(Z_{\ell})},\quad \sigma_{+} = v_{r} + \sqrt{\frac{\vrs_{r}}{\vr_{r} (\vrs_{r} - \vr_{r})}}\sqrt{\bar p -p_{0}(Z_{r})}.
\end{equation*}
\end{prop}
This proposition can be proven using similar arguments as in \cite{DeNaBonSan}.

\begin{figure}
\includegraphics[width=0.5\textwidth]{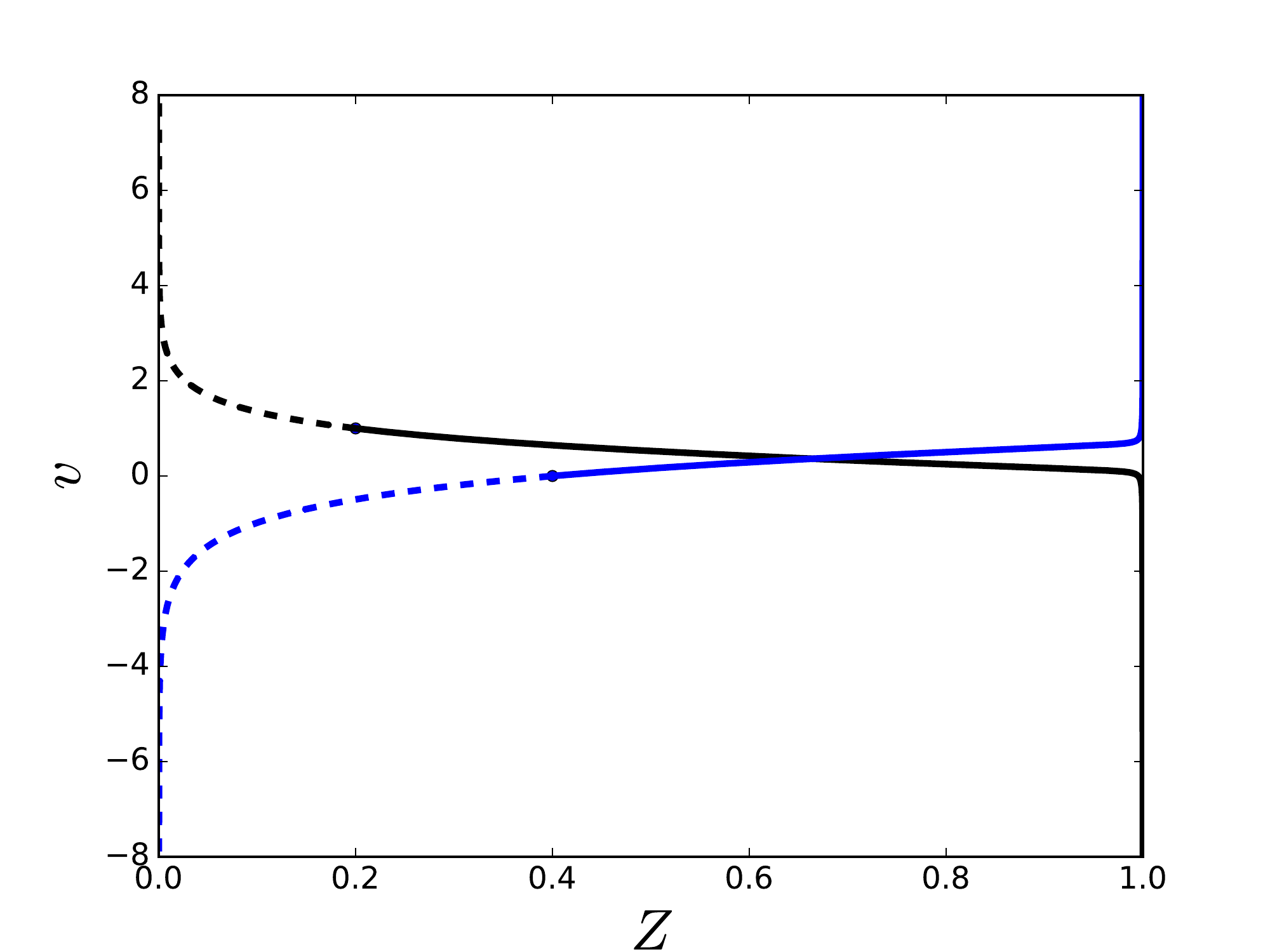}
\includegraphics[width=0.5\textwidth]{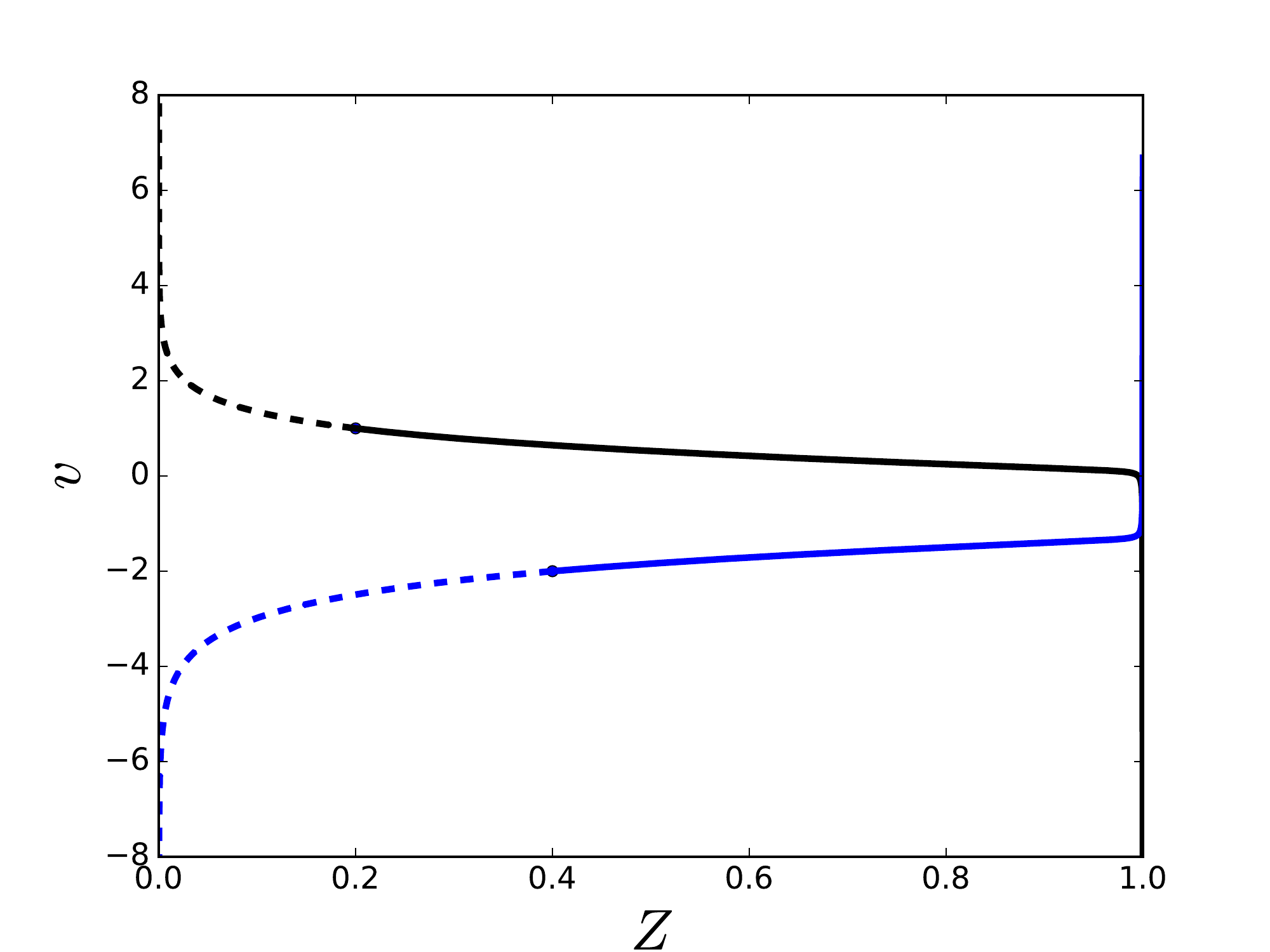}
\caption{Intersection of the 1-integral/Hugoniot curve issued from the left state $(\vr_{\ell}, v_{\ell}, Z_{\ell}) = (0.8, 1, 0.2)$ and the 3-integral/Hugoniot curve issued from the right state for $\eps = 10^{-3}$. The rarefaction curves are in dashed line and  the shock curve in solid line. Left: the right state is given by $(\vr_{r}, v_{r}, Z_{r}) = (0.8, 0, 0.4)$ and the intermediate state $(v_{m}^{0}, Z_{m}^{0})$ is not a congested state. Right: the right state is given by $(\vr_{r}, v_{r}, Z_{r}) = (0.8, -2, 0.4)$ and the intersection point is very closed to the congested line $Z = 1$.}
\label{fig:Riem}
\end{figure}

Below, on Figure \ref{fig:Riem} we present two different solutions to the Riemann problem \eqref{eq:system1D}-\eqref{data1D}. Depending on the initial location of the left and right states, the intersection state $(v_m,Z_m)$ might be a congested state or not.

\section{Fully discrete scheme in dimension 2}
\label{sec:scheme_2d}

We consider the computational domain $[0,1]\times [0,1]$ and spatial space steps $\Delta x = 1/N_{x}, \Delta y = 1/N_{y} > 0$, with $N_{x}, N_{y} \in \mathbb{N}$: the mesh points are thus $\boldsymbol{x}_{i,j} = (i \Delta x, j\Delta y)$, $\forall (i,j) \in \left\{0,\ldots, N_{x}\right\}\times\left\{0,\ldots, N_{y}\right\}$. Let $\vr^{n}_{i,j}$, $\q^{n}_{i,j}$, $Z^{n}_{i,j}$, $\vr^{\ast\, n}_{i,j}$ denote the approximate solution at time $t^{n}$ on mesh cell $[i \Delta x, (i+1)\Delta x]\times [j \Delta x, (j+1)\Delta x]$. 

The two-dimensional version of \eqref{eq:fulldiscrete} reads:
\begin{align}
&\frac{\vr^{n+1}_{i,j} - \vr^{n}_{i,j}}{\Delta t} + \frac{1}{\Delta x} (F_{(i+\frac{1}{2},j)}^{n+1} - F_{(i-\frac{1}{2},j)}^{n+1}) + \frac{1}{\Delta y}  (\widetilde F_{(i,j+\frac{1}{2})}^{n+1} - \widetilde F_{1,(i,j-\frac{1}{2})}^{n+1}) = 0,\\
&\frac{\q^{n+1}_{i,j} - \q^{n}_{i,j}}{\Delta t} + \frac{1}{\Delta x} (\boldsymbol{G}_{(i+\frac{1}{2},j)}^{n} - \boldsymbol{G}_{(i-\frac{1}{2},j)}^{n}) + \frac{1}{\Delta y}  (\widetilde{\boldsymbol{ G}}_{(i,j+\frac{1}{2})}^{n} -  \widetilde{\boldsymbol{ G}}_{(i,j-\frac{1}{2})}^{n})\nonumber\\
&\hspace{7cm} + (\nabla\pi_\ep(Z^{n+1}))_{i,j} = 0,\label{eq:discrete_q_2d}\\
&\frac{Z^{n+1}_{i,j} - Z^{n}_{i,j}}{\Delta t} + \frac{1}{\Delta x} (H_{(i+\frac{1}{2},j)}^{n+1} - H_{(i-\frac{1}{2},j)}^{n+1})+ \frac{1}{\Delta y}  (\widetilde H_{(i,j+\frac{1}{2})}^{n+1} - \widetilde H_{3,(i,j-\frac{1}{2})}^{n+1})  = 0\label{eq:discrete_z_2d}.
\end{align}
where fluxes $F^{n+1}$, $\boldsymbol{G}^{n}$, $H^{n+1}$ (in the first spatial direction) are defined:
\begin{align}
&F_{(i+\frac{1}{2},j)}^{n+1} = \frac{1}{2}\big(q^{n+1}_{1,(i+1,j)}+q^{n+1}_{1,(i,j)}\big) - (D_{\vr})_{i+\frac{1}{2},j}^{n},\\
&\boldsymbol{G}_{(i+\frac{1}{2},j)}^{n} = \frac{1}{2}\big(\boldsymbol{f}_{(i+1,j)}^{n}+\boldsymbol{f}^{n}_{(i,j)}\big) - (\boldsymbol{D}_{\q})_{i+\frac{1}{2},j}^{n},\\
&H_{(i+\frac{1}{2},j)}^{n+1} = \frac{1}{2}\Big( \frac{Z^{n}_{i+1,j}}{\vr^{n}_{i+1,j}} q^{n+1}_{1,(i+1,j)}+ \frac{Z^{n}_{i,j}}{\vr^{n}_{i,j}} q^{n+1}_{1,(i,j)}\Big) - (D_{Z})_{i+\frac{1}{2},j}^{n},,\label{eq:2space_H}
\end{align}
with \begin{equation*}
\boldsymbol{f}^{n} = \begin{bmatrix} (q_{1}^{n})^{2} + p(Z^{n})\\
q_{1}^{n} q_{2}^{n}
\end{bmatrix}. 
\end{equation*}
Fluxes $\widetilde F^{n+1}$, $\widetilde{\boldsymbol{G}}^{n}$, $\widetilde{H}^{n+1}$ in the second spatial direction are defined by:
\begin{align}
&\widetilde F_{(i,j+\frac{1}{2})}^{n+1} = \frac{1}{2}\big(q^{n+1}_{2,(i,j+1)}+q^{n+1}_{2,(i,j)}\big) - (D_{\vr})_{i,j+\frac{1}{2}}^{n},\\
&\widetilde{\boldsymbol{G}}_{(i,j+\frac{1}{2})}^{n} = \frac{1}{2}\big(\widetilde{\boldsymbol{f}}_{(i,j+1)}^{n}+\widetilde{\boldsymbol{f}}^{n}_{(i,j)}\big)- (\boldsymbol{D}_{\q})_{i,j+\frac{1}{2}}^{n},\\
&\widetilde H_{(i,j+\frac{1}{2})}^{n+1} = \frac{1}{2}\Big( \frac{Z^{n}_{i,j+1}}{\vr^{n}_{i,j+1}} q^{n+1}_{2,(i,j+1)}+ \frac{Z^{n}_{i,j}}{\vr^{n}_{i,j}} q^{n+1}_{2,(i,j)}\Big) - (D_{Z})_{i,j+\frac{1}{2}},\label{eq:2space_tildeH}
\end{align}
with \begin{equation*}
\widetilde{\boldsymbol{f}}^{n} = \begin{bmatrix} 
q_{1}^{n} q_{2}^{n}\\
(q_{2}^{n})^{2} + p(Z^{n})\\
\end{bmatrix}. 
\end{equation*}
The upwindings $D_{\vr}$, $\boldsymbol{D}_{\q}$, $D_{Z}$ are defined similarly as for the one-dimensional case (sse \eqref{eq:upwindings}-\eqref{eq:maxcharspeed}).

The implicit pressure in \eqref{eq:discrete_q_2d}  is discretized by the centered difference:
\begin{equation*}
(\nabla\pi_\ep(Z^{n+1}))_{i,j} = \begin{bmatrix}
\displaystyle\frac{\pi_\ep(Z^{n+1}_{i+1,j}) - \pi_\ep(Z^{n+1}_{i-1,j})}{2\Delta x}\\
\displaystyle\frac{\pi_\ep(Z^{n+1}_{i,j+1}) - \pi_\ep(Z^{n+1}_{i,j-1})}{2\Delta y}
\end{bmatrix}.
\end{equation*}
Inserting equation \eqref{eq:discrete_q_2d} into \eqref{eq:discrete_z_2d}, we obtain:
\begin{align*}
&Z^{n+1}_{i,j} - Z^{n}_{i,j} + \frac{\Delta t}{\Delta x} \big(\bar H_{(i+\frac{1}{2},j)}^{n} - \bar H_{(i+\frac{1}{2},j)}^{n}\big) + \frac{\Delta t}{\Delta y} \big( \bar{\widetilde H}_{(i+\frac{1}{2},j)}^{n} -  \bar{\widetilde H}_{(i+\frac{1}{2},j)}^{n}\big) \\
&- \frac{\Delta t^{2}}{\Delta x^{2}} \frac{1}{2}\Big( \frac{Z^{n}_{i+1,j}}{\vr^{n}_{i+1,j}} \big(G_{(i+\frac{3}{2},j),1}^{n} - G_{(i+\frac{1}{2},j),1}^{n}\big)- \frac{Z^{n}_{i-1,j}}{\vr^{n}_{i-1,j}} \big(G_{(i-\frac{1}{2},j),1}^{n} - G_{(i-\frac{3}{2},j),1}^{n}\big)\Big)\\
&- \frac{\Delta t^{2}}{\Delta x \Delta y} \frac{1}{2}\Big( \frac{Z^{n}_{i+1,j}}{\vr^{n}_{i+1,j}} \big(\widetilde G_{(i+1,j+\frac{1}{2}),1}^{n} - \widetilde G_{(i+1,j-\frac{1}{2}),1}^{n}\big) \\
&\hspace{5cm}	- \frac{Z^{n}_{i-1,j}}{\vr^{n}_{i-1,j}} \big(\widetilde G_{(i-1,j+\frac{1}{2}),1}^{n} - \widetilde G_{(i-1,j-\frac{1}{2}),1}^{n}\big)\Big)\\
&-\frac{\Delta t^{2}}{\Delta y^{2}} \frac{1}{2}\Big( \frac{Z^{n}_{i,j+1}}{\vr^{n}_{i,j+1}} \big(\widetilde G_{(i,j+\frac{3}{2}),2}^{n} - \widetilde G_{(i,j+\frac{1}{2}),2}^{n}\big)- \frac{Z^{n}_{i,j-1}}{\vr^{n}_{i,j-1}} \big(\widetilde G_{(i,j-\frac{1}{2}),2}^{n} - \widetilde G_{(i,j-\frac{3}{2}),2}^{n}\big)\Big)\\
&- \frac{\Delta t^{2}}{\Delta x \Delta y} \frac{1}{2}\Big( \frac{Z^{n}_{i,j+1}}{\vr^{n}_{i,j+1}} \big(G_{(i+\frac{1}{2},j+1),2}^{n} - G_{(i-\frac{1}{2},j+1),2}^{n}\big) \\
&\hspace{5cm}	- \frac{Z^{n}_{i,j-1}}{\vr^{n}_{i,j-1}} \big( G_{(i+\frac{1}{2},j-1),2}^{n} -  G_{(i-\frac{1}{2},j-1),2}^{n}\big)\Big)\\
&-  \frac{\Delta t^{2}}{\Delta x^{2}} \frac{1}{4}\Big( \frac{Z^{n}_{i+1,j}}{\vr^{n}_{i+1,j}} \big(\pi_\ep(Z^{n+1}_{i+2,j}) - \pi_\ep(Z^{n+1}_{i,j})\big) - \frac{Z^{n}_{i-1,j}}{\vr^{n}_{i-1,j}} \big(\pi_\ep(Z^{n+1}_{i,j}) - \pi_\ep(Z^{n+1}_{i-2,j})\big)\Big) \\
&-  \frac{\Delta t^{2}}{\Delta y^{2}} \frac{1}{4}\Big( \frac{Z^{n}_{i,j+1}}{\vr^{n}_{i,j+1}} \big(\pi_\ep(Z^{n+1}_{i,j+2}) - \pi_\ep(Z^{n+1}_{i,j})\big) - \frac{Z^{n}_{i,j-1}}{\vr^{n}_{i,j-1}} \big(\pi_\ep(Z^{n+1}_{i,j}) - \pi_\ep(Z^{n+1}_{i,j-2})\big)\Big) = 0,
\end{align*}
where terms $\bar H^{n}$ and $ \bar{\widetilde H}^{n}$ have the same expressions as \eqref{eq:2space_H}-\eqref{eq:2space_tildeH} but where all quantities are taken
explicitly.

 \section{Second order in time $(\vr,q)$-method/SL}
\label{sec:secondorder_rhoscheme}

The second order accuracy scheme for the $(\vr,q)$-method/SL is based on a Strang splitting between advection of congestion density and advection of $(\vr, \q)$. It consists in the following steps:
\begin{enumerate}
\item Compute $\vr^{\ast\ n+1/2}$ by solving the advection of over $\Delta t/2$
\begin{align*}
\frac{{\vr^*}^{n+1/2} - {\vr^*}^{n}}{\Delta t/2} + \frac{\q^{n}}{\vr^{n}}\cdot\nabla_{x} {\vr^*}^{n} = 0.
\end{align*}
\item Compute $(\vr^{n+1}, q^{n+1})$ with the RK2CN scheme as proposed in \cite{CorDeKum}:

\textbf{First step (half time step):}
\begin{subequations}
\begin{align*}
&\frac{\vr^{n+1/2} - \vr^{n}}{\Delta t/2} + \nabla_{x} \cdot \q^{n+1/2} - \mathcal{D}_{\vr}^{n} = 0,\\
&\frac{\q^{n+1/2} - \q^{n}}{\Delta t/2} + \nabla_{x}\cdot\left(\frac{\q^{n}\otimes \q^{n}}{\vr^{n}} + p(Z^{n}) \Id \right) -  \mathcal{D}_{q}^{n} + \nabla_{x}(\pi_\ep(\vr^{n+1/2}/\vr^{\ast, n+1/2})) = 0.
\end{align*}
\end{subequations}  
\textbf{Second step (full time step):}
\begin{subequations}
\begin{align*}
&\frac{\vr^{n+1} - \vr^{n}}{\Delta t} + \nabla_{x} \cdot \left(\frac{\q^{n+1}+\q^{n}}{2}\right) - \mathcal{D}_{\vr}^{n} = 0,\\
&\frac{\q^{n+1} - \q^{n}}{\Delta t} + \nabla_{x}\cdot\left(\frac{\q^{n+1/2}\otimes \q^{n+1/2}}{\vr^{n+1/2}} + p(Z^{n+1/2}) \Id \right) -  \mathcal{D}_{q}^{n+1/2} \nonumber\\
&\hspace{3cm}+ \nabla_{x}\left(\frac{\pi_\ep(\vr^{n}/\vr^{\ast, n+1/2}) +\pi_\ep(\vr^{n+1}/\vr^{\ast, n+1/2})}{2}\right) = 0.
\end{align*}
\end{subequations}  

where $\mathcal{D}_{\vr}$, $\mathcal{D}_{q}$ denote the numerical diffusion coming from fluxes.
\item Advection of $\vr^{\ast}$ on $\Delta t /2$ time step
\begin{align*}
\frac{{\vr^*}^{n+1} - {\vr^*}^{n+1/2}}{\Delta t/2} + \frac{\q^{n+1}}{\vr^{n+1}}\cdot\nabla_{x}  {\vr^*}^{n+1/2} = 0.
\end{align*}
\end{enumerate}
A second order in time version of the semi-Lagrangian scheme has to be used. We here consider the second order Taylor approximation of the caracteristic line whose one-dimensional version reads:
\begin{align*}
\vr_{i}^{\ast n+1} = \left[\Pi \vr^{\ast n}\right]\Big(x_{i} - v_{i} \Delta t +  a_{i} v_{i} \frac{\Delta t^{2}}{2} \Big).
\end{align*}
where $v_{i} =  q_{i}/\vr_{i}$ for all $i$ and $a_{i}$ is an upwind finite difference approximation of the first derivative of the velocity: $a_{i} = (v_{i} - v_{i-1})/\Delta x$ if $v_{i} > 0$ and $a_{i} = (v_{i+1} - v_{i})/\Delta x$ if $v_{i} \leqslant 0$. 

\bigskip

\section*{Acknowledgements}
P.D. acknowledges support by the Engineering and Physical Sciences Research Council (EPSRC) under grant no. EP/M006883/1, by the Royal Society and the Wolfson Foundation through a Royal Society Wolfson Research Merit Award no. WM130048 and by the National Science Foundation (NSF) under grant no. RNMS11-07444 (KI-Net). P.D. is on leave from CNRS, Institut de Math\'ematiques de Toulouse, France. P.M. acknowledges the support of MAThematics Center Heidelberg (MATCH). E.Z. was supported by the the Department of Mathematics, Imperial College, through a Chapman Fellowship, and by the Polish Ministry of Science and Higher Education grant Iuventus Plus  no. 0888/IP3/2016/74. 
\medskip

\section*{Data availability}
No new data were collected in the course of this research.

\end{document}